%% file: main.tex
\documentclass{article}

\usepackage{epsfig,amsmath,amssymb,amsfonts}
\usepackage{cite,placeins}

\newtheorem{Definition}{\textbf{Definition}}[section]
\newtheorem{Theorem}{\textbf{Theorem}}[section]
\newtheorem{Conjecture}{\textbf{Conjecture}}[section]
\newtheorem{Lemma}{\hskip 0pt Lemma}[section]
\newtheorem{Corollary}{\hskip 0pt Corollary}[section]
\newtheorem{Proposition}{\hskip 0pt Proposition}[section]
\newenvironment{Proof}{\medskip\par
{{\textit{\textbf{Proof:}}}\,}}{{\mbox{\,$\Box$}\par}}

\title{On embedded spheres of affine manifolds}
\author{Weiqiang Wu\\
Department of Mathematics\\
University of Maryland, College Park\\
\texttt{waikng@math.umd.edu}}

\begin{document}
\maketitle

\begin{abstract}
This paper studies certain embedded spheres in closed affine manifolds. For $n \geq 3$, we investigate the dome bodies in a closed affine $n$-manifold $M$ with its boundary homeomorphic to a sphere under the assumption that a developing map restricted to a component of $\partial\hat{M}$ is an embedding onto a strictly convex sphere in $\mathbb{A}^n$. By using the recurrent property of an incomplete geodesic we show that dome bodies are compact. Then a maximal dome body is a closed solid ball bounded by a component of $\partial\hat{M}$, and hence equals $\hat{M}$. The main theorem is that the standard ball in an affine space can only bound one compact affine manifold inside, namely the solid ball.
\end{abstract}

\include{Chapter1}
\include{Chapter2}
\include{Chapter3}

\include{Chapter4}

\small\normalsize
\bibliographystyle{amsalpha}
\bibliography{affine}

\end{document}

%% file: Chapter1.tex

\section{Introduction}
Following Felix Klein's 1872 Erlanger program, geometry is the study of properties of a space $X$ invariant under a group $G$ of transformations of $X$ \cite{Projective}. Classical geometries include Euclidean geometry $\left(\mathbb{E}^n, E(n)\right)$, spherical geometry $\left(\mathbb{S}^n, O(n+1)\right)$ and hyperbolic geometry $\left(\mathbb{H}^n, PO(n,1)\right)$ \cite{Foundations}. The transformation groups of these geometries preserve metrics with constant curvature $0$, $1$ and $-1$ respectively, so in particular the techniques in Riemannian geometry can be applied. Projective
geometry $\left(\mathbb{RP}^n, PGL(n+1)\right)$ unifies all three geometries into a more general category. Though no longer a metrical form of geometry, it greatly benefits from the classical geometries, and especially it is largely motivated by hyperbolic geometry. The affine geometry $\left(\mathbb{A}^n, \mbox{Aff}(n) \right)$, as the intermediate layer between Euclidean geometry and projective geometry, on the contrary, is not easy to study.

\begin{figure}[h]
\begin{center}
\scalebox{0.5} {\epsffile{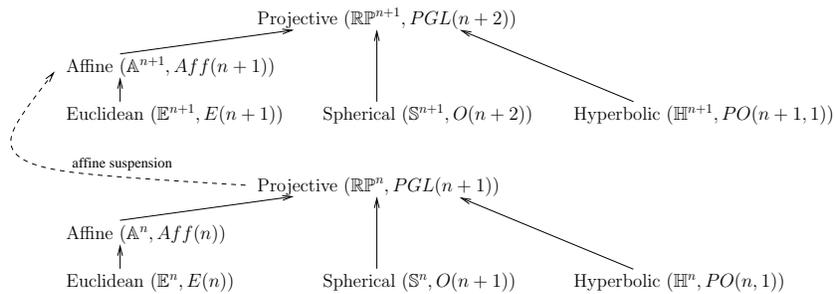}}
\end{center}
\caption{Relations between the geometries.} \label{fig:xg_hier}
\end{figure}

\begin{figure}[h]
\begin{center}
\scalebox{0.6} {\epsffile{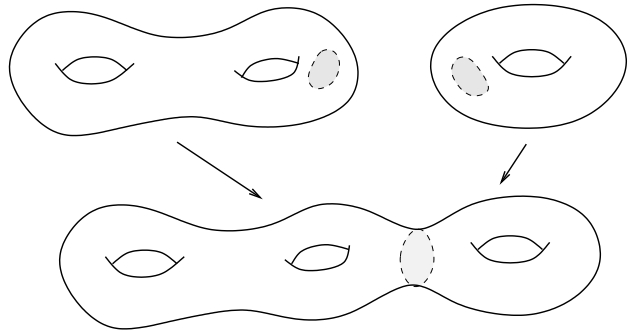}}
\end{center}
\caption{A connected sum.} \label{fig:connected_sum}
\end{figure}

This paper studies certain embedded spheres of closed affine manifolds with intent to understand the topological structures of affine manifolds. More specifically, we are interested in the following basic topological question: \begin{center}\it ``Can a non-trivial connected sum of closed manifolds admit an affine structure?''\end{center}

This question is somewhat related to the famous and infamous
\begin{Conjecture}[Chern conjecture]
The Euler class of a closed affine manifold vanishes.
\end{Conjecture}
Kostant and Sullivan proved it in the complete case \cite{Kostant-Sullivan_1975}. Several special cases were also worked out by others; however it is notoriously hard in general.

Smillie constructed some interesting examples of flat manifolds with non-zero Euler characteristic from connected sum \cite{Smillie_1977}, which might potentially be counter-examples to the {\it Chern conjecture}.

Instead of attacking the original question, our first attempt is to answer a simpler one: \begin{center}\it ``Can the standard sphere bound a compact affine manifold inside other than the solid ball?''\end{center}

\begin{figure}[h]
\begin{center}
\scalebox{0.6} {\epsffile{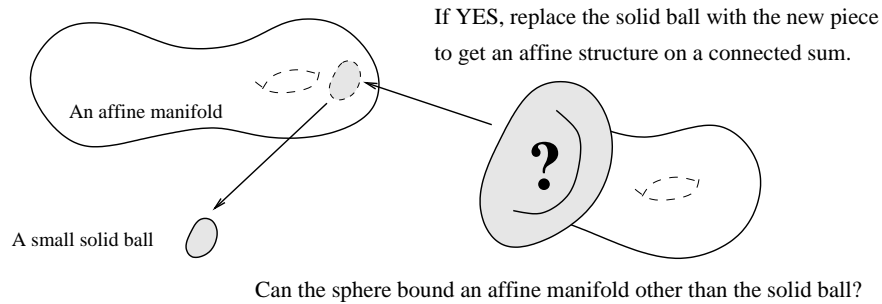}}
\end{center}
\caption{Motivation for the simpler question.} \label{fig:S_bounds_a_ball}
\end{figure}

We adopt similar geometric techniques as in \cite{Fried_1980, Carriere_1989, Choi_decomp_1999}. By looking at some geometric objects together with the recurrence of an incomplete geodesic, we can give a negative answer to the above question. More specifically we prove
\begin{Theorem}\label{thm:main}
    For $n \geq 3$, let $(M,\partial)$ be a compact affine $n$-manifold with boundary $\partial \simeq S^{n-1}$ and $dev: \hat{M} \to \mathbb{A}^n$ be a developing map. If
    \begin{itemize}
    \item $dev$ restricted to some lift $\hat{\partial}$ of $\partial$ is an embedding,
    \item $dev$ maps a neighborhood of $\hat{\partial}$ to the closure of the bounded part of $\mathbb{A}^n \backslash dev(\hat{\partial})$,
    \end{itemize}
    then $(M, \partial)$ is homeomorphic to $(D^n, S^{n-1})$.
\end{Theorem}

Then the following corollaries concerning embedded spheres of closed affine manifolds follow directly.
\begin{Corollary}\label{cor:S_bounds_B}
    For $n \geq 3$, let $M$ be a closed affine $n$-manifold and $S$ be an embedded separating $(n$-$1)$-sphere. If a developing map $dev$ restricted to some lift $\hat{S}$ of $S$ is an embedding, then $S$ bounds a $n$-ball in $M$.
\end{Corollary}

\begin{Corollary}\label{cor:S_non_injective}
    For $n \geq 3$, a developing map restricted to a lift of a non-trivial separating sphere of a closed affine manifold is not injective.
\end{Corollary}

Though the theme of this paper is affine structures, we can extend Theorem \ref{thm:main} to the projective case.
\begin{Theorem}\label{thm:proj}
    For $n \geq 3$, let $(M, \partial)$ be a compact projective n-manifold with boundary $\partial$ homeomorphic to $S^{n-1}$. If
    \begin{itemize}
      \item $dev$ restricted to some lift $\hat{\partial}$ of $\partial$ is an embedding,
      \item $dev(\hat{\partial})$ is contained in an affine patch $\mathbb{A}^n$,
      \item $dev$ maps a neighborhood of $\hat{\partial}$ to the closure of the bounded part of $\mathbb{A}^n \backslash dev(\hat{\partial})$,
    \end{itemize}
    then $(M, \partial)$ is homeomorphic to $(D^n, S^{n-1})$.
\end{Theorem}

 This paper is organized as follow:

 Section 2 gives a quick review of the general theory of $(X,G)$-manifolds. We introduce the pair of developing map and holonomy, which plays as an important tool to study $(X,G)$-manifolds in general. Then we list how examples of $(X,G)$-manifolds arise. After that we briefly review the three classical geometries: Euclidean geometry, spherical geometry and hyperbolic geometry, and provide examples of these geometries in dimensions $\leq$ 3.

 Section 3 focuses on affine structures and provides closed orientable examples in dimension 2 and 3. We summarize the results of the classification of affine structures on the 2-torus, the classification of complete affine 3-manifolds and the classification of radiant affine 3-manifolds. We also discuss projective structures and the construction of affine suspension from projective manifolds in the course. At the end we discuss Goldman's two non-complete and non-radiant examples. Basically, we provide all the known examples in dimension 3.
 
 Section 4 uses the generalized Schoenflies theorem to reduce Theorem \ref{thm:main} to Theorem \ref{thm:convex_main}. Then we develop the theory of dome bodies to prove Theorem \ref{thm:convex_main}. We adopt similar techniques used in \cite{Fried_1980, Carriere_1989, Choi_decomp_1999} to reduce the proof to the main technical point: the compactness of dome bodies. Lastly, we extend the proof to the projective case. 

%% file: Chapter2.tex

\section{(X,G)-manifolds}\label{Ch2}

In this paper we work only in the $C^{\infty}$ category. Without further specification manifolds and maps are assumed to be smooth.

\subsection{Group actions}
We recall few definitions from the theory of group actions in this subsection for the reader's convenience.

\begin{Definition}[$G$ acts on $X$]
Let $X$ be a manifold. A Lie group $G$ acts on $X$ via $\Phi$, if $\Phi: G \times X \to X$ is a smooth map such that
$$\Phi(1, x) = x \quad \mbox{and} \quad \Phi(g, \Phi(h, x)) = \Phi(g.h, x).$$
for all $g, h \in G$ and $x \in X$. In other words, $G$ acts on $X$ via $\Phi$ if $\Phi$ induces an homomorphism from $G$ to $\mbox{Diff}(X)$.
\end{Definition}

Remark. We usually omit $\Phi$ if $G$ acts on $X$ in some natural way. We also use $g.x$ for $\Phi(g,x)$ and $G.x$ for the orbit of $x$ ( i.e. $G.x = \left\{\, g.x \,\middle|\, g \in G \,\right\}$ ).

\begin{Definition}[Free action]
$G$ acts on $X$ freely, if $g.x \neq x$ for all $g \in G-\{1\}$ and $x \in X$.
\end{Definition}

\begin{Definition}[Transitive action]
$G$ acts on $X$ transitively, if there exists some ( and hence for all ) $x \in X$, such that $G.x = X$.
\end{Definition}

\begin{Definition}[Proper action]
$G$ acts on $X$ properly, if the map $$\Phi \times I: G \times X \to X \times X \quad \quad (g, x) \mapsto (g.x, x)$$
is proper, i.e. for any compact set $K$ of $X \times X$, the inverse image $(\Phi \times I)^{-1}(K)$ is compact.
\end{Definition}

\begin{Definition}[Properly discontinuous action]
A discrete group $\Gamma$ acts on $X$ properly discontinuously, if for any compact set $K \subset X$, the set $$\left\{\, \gamma \in \Gamma \,\middle|\, \gamma.K \cap K \neq \emptyset \,\right\}$$ is finite.
\end{Definition}

Remark. That $G$ acts properly discontinuously on $X$ is the same as that $G$ with the discrete topology acts properly on $X$.

\newpage
\subsection{(X,G)-manifolds}
Let $X$ be a connected manifold, and $G$ acts transitively on $X$. Given any manifold $M$ of the same dimension as $X$, we want to give $M$ an $(X,G)$-structure as follows:
\begin{Definition}[(X,G)-structure]
An $(X,G)$-structure on $M$ is given by an open cover $\{U_{\alpha}\}$ of $M$ together with the charts $\{ \phi_{\alpha}: U_{\alpha} \to X\}$ which are homeomorphisms onto their images, and there exists $\{ g_{\alpha \beta} \in G\}$ such that $$g_{\alpha \beta}\circ \phi_{\beta} = \phi_{\alpha} \quad \mbox{on} \quad  U_{\alpha}\cap U_{\beta}$$
$$g_{\alpha \beta}.g_{\beta \gamma} = g_{\alpha \gamma} \quad \mbox{and} \quad g_{\alpha \alpha} = 1.$$
$M$ is said to be an $(X,G)$-manifold, if it is given an $(X,G)$-structure.
\end{Definition}

\begin{figure}[h]
\begin{center}
\scalebox{0.5} {\epsffile{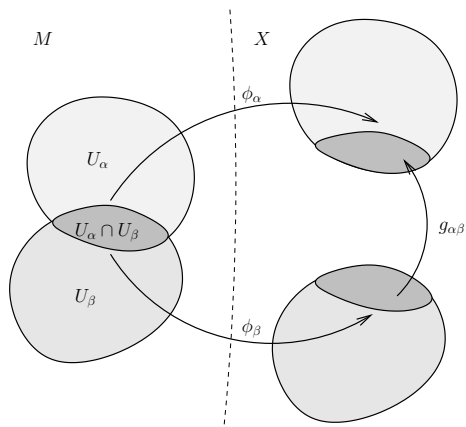}}
\end{center}
\caption{Charts for an (X,G)-structure.} \label{fig:Charts}
\end{figure}

Remark. The geometry of interest is $G$-invariant and different charts in an $(X,G)$-structure are related by a $g$ action. Then an $(X,G)$-structure on $M$ pulls the local geometry of $X$ to $M$ via the charts.

Note that if we pull back the charts by a covering map, we can get an $(X,G)$-structure on the covering space. Hence we have
\begin{Proposition}[(X,G)-structure can be lifted]
Let $M' \to M$ be a covering map. If $M$ is given an $(X,G)$-structure, then there is a compatible $(X,G)$-structure on $M'$. \label{prop:lift_structure}
\end{Proposition}

\subsubsection{Developing map and holonomy}
An important and powerful tool is the pair consisting of a developing map and the corresponding holonomy representation. Before defining them, we need the following property for the $G$ action.
\begin{Definition}[Strongly effective action]
$G$ acts on $X$ strongly effectively, if whenever $g$ fixes $U$ pointwise for some open set $U \subset X$, then $g=1$.
\end{Definition}

Remark. All the $(X,G)$-geometries in this paper have this property. It is also worth noting that $\mbox{Diff}(X)$ acting on $X$ does not have this property.

Given an $(X,G)$-structure on $M$, the universal cover $\hat{M}$ admits an $(X,G)$-structure by Proposition \ref{prop:lift_structure}. We then take a base point in $\hat{M}$ and an $(X,G)$-chart containing the base point and mimic the construction of analytic continuation in complex analysis. The strong effectiveness guarantees that $dev$ is well defined, i.e. independent of the path chosen to do the continuation.

\begin{figure}[h]
\begin{center}
\scalebox{0.4} {\epsffile{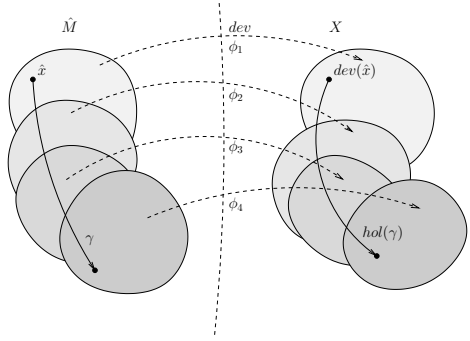}}
\end{center}
\caption{Developing map and holonomy.} \label{fig:Dev}
\end{figure}

\begin{Definition}[Developing map and holonomy]
Given an $(X,G)$-structure on $M$, there exists a pair $(dev, hol)$,
$$ dev: \hat{M} \to X, \quad \quad hol: \pi_1(M) \to G$$
$$dev(\gamma.\hat{x}) = hol(\gamma).dev(\hat{x})$$ for all $\gamma \in \pi_1(M)$ and $\hat{x} \in \hat{M}$. $dev$ is a local homeomorphism and $hol$ is a homomorphism. The pair is unique up to the conjugation of an element $g \in G$.
\end{Definition}

Remark. A developing pair $(dev, hol)$ actually determines the $(X,G)$-structure.

The usefulness of development is readily seen in the following simple proposition.
\begin{Proposition}
Let $M$ be a manifold with compact universal cover. If $X$ is not compact, then $M$ admits no $(X,G)$-structure.
\end{Proposition}
\begin{Proof}
Suppose $M$ admits an $(X,G)$-structure. Then $dev(\hat{M})$ is both open and compact, and hence equals $X$. This contradicts to that $X$ is not compact.
\end{Proof}

\subsubsection{Examples of (X,G)-manifolds}\label{subsec:examples}
We list few typical examples of $(X,G)$-manifolds in the following:
\begin{itemize}
\item The fundamental example: $M = X$.

 $X$ itself together with the identity map as a chart is a natural $(X,G)$-structure.

\item The space forms: $M = X/\Gamma$.

 Here $\Gamma \subset G$ is a discrete subgroup that acts freely and properly discontinuously on $X$. The set of local inverses of the covering charts $X \to X/\Gamma$ gives a natural $(X,G)$-structure.

 Remark. The reason we are interested in a discrete subgroup $\Gamma$ acting freely and properly discontinuously on $X$ is that the quotient space $X/\Gamma$ then has a manifold structure automatically. If $\Gamma$ further acts cocompactly, the quotient space is a closed manifold.

 In this case, $dev: \hat{M} \to X$ is a covering map. If $X$ is simply connected, then $\hat{M}$ is homeomorphic to $X$ and $hol: \pi_1(M) \to \Gamma$ is an isomorphism. We define this as

 \begin{Definition}[Complete structure]
   An $(X,G)$-structure on $M$ is called complete, if $dev: \hat{M} \to X$ is a covering map. \label{def:Complete}
 \end{Definition}

 Remark. In the case that $X$ has a $G$-invariant Riemannian metric, the completeness defined above of the $(X,G)$-structure on $M$ is equivalent to the completeness of $M$ as a metric space induced from the Riemannian metric on $X$. Therefore all three classical geometric structures ( Euclidean, projective and hyperbolic ) on a closed manifold belong to this type by the Hopf-Rinow theorem.

\item $M = \hat{U}/\Gamma$, where $U$ is a connected proper domain in $X$ and $\hat{U}$ is the universal cover of $U$.

 Similar to the previous case, $\Gamma \subset G$ is a discrete subgroup that acts freely and properly discontinuously on $\hat{U}$. The set of the inverses of the covering charts $\hat{U} \to  \hat{U}/\Gamma$ gives a natural $(X,G)$-structure. In this case $dev: \hat{M} \to X$ is a covering onto $U$.

 This type of example arises in affine and projective geometries. An example arises by considering a closed hyperbolic surface as a projective surface. They provide building blocks for projective structures on surfaces.

\item $dev: \hat{M} \to X$ is not a covering onto its image.

 This type of example also arises in affine and projective geometries. An example of this type will be seen in \S 3.3.2.2.
\end{itemize}

\subsection{Euclidean manifolds}
\begin{Definition}[Euclidean structure]
$(X,G)$-structures are called Euclidean, when
$$ X = \mathbb{E}^n = \left\{\, (x_0, \ldots, x_{n-1})^T \,\middle|\, x_i \in \mathbb{R} \,\right\}$$
$$ G = E(n) = \mathbb{R}^n \rtimes O(n),$$
where $G$ acts on $X$ by
$$ g.x = A \cdot x + b$$
for $g = (b, A) \in G$ and $x \in X$.
\end{Definition}

Remark. $E(n)$ preserves a Riemannian metric on $\mathbb{E}^n$
$$ds^2 = dx_0^2 + dx_1^2 + \ldots + dx_{n-1}^2$$
with curvature $0$. Hence Euclidean manifolds has Euler characteristics $0$ by Gauss-Bonnet-Chern theorem.

What makes Euclidean manifolds easy to classify is that the linear part of $E(n)$ is the compact group $O(n)$ and discrete subgroups of $E(n)$ have some nice properties:

\begin{Lemma}[\S5.4 Lemma 4 and 5 in \cite{Foundations}]
Let $\Gamma$ be a discrete subgroup of $E(n)$ and let $\phi = a + A$ and $\psi = b + B$ be in $\Gamma$ with $\| A-I \| < 1$ and $\|B-I \| < 1$. Then $\phi$ and $\psi$ commute.
\end{Lemma}

\begin{Lemma}[Theorem 5.4.6 in \cite{Foundations}]
Let $\Gamma$ be a discrete subgroup of $E(n)$. Then $\Gamma$ has a free abelian subgroup $H$ of rank $m$ and finite index.
\end{Lemma}

Applying the above lemma, we see
\begin{Theorem}
Every closed Euclidean manifold is finitely covered by a torus.
\end{Theorem}

Furthermore we have
\begin{Theorem}[Bieberbach's theorem, Theorem 7.5.3 in \cite{Foundations}] \label{thm:Bieberbach}
There are only finitely many n-dimensional Euclidean manifolds up to affine equivalence for each $n$.
\end{Theorem}

\subsubsection{Closed orientable examples in low dimensions}
\begin{figure}[h]
\begin{center}
\scalebox{0.7} {\epsffile{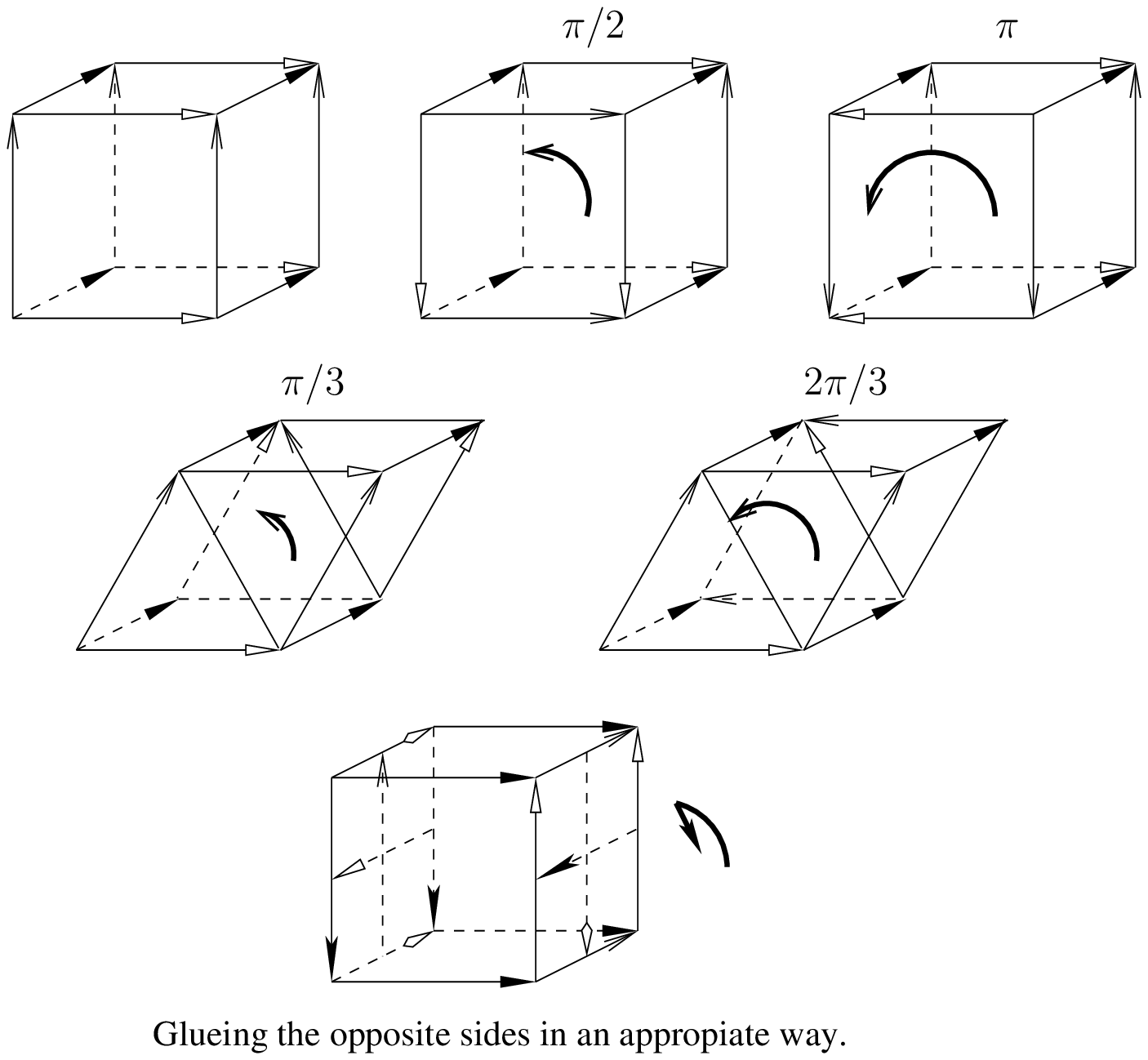}}
\end{center}
\caption{The 6 affine equivalence classes of closed orientable Euclidean manifolds.} \label{fig:euclid3}
\end{figure}

\begin{itemize}
\item Dim = 1
    \begin{itemize}
    \item $T^1 = E^1/\Gamma$, where $\Gamma$ is cyclic subgroup generated by a translation.
    \end{itemize}
\item Dim = 2
    \begin{itemize}
    \item $T^2 = E^2/\Gamma$, where $\Gamma$ is an abelian subgroup of rank $2$ generated by two independent translations.
    \end{itemize}
\item Dim = 3
    \begin{itemize}
    \item $T^3 / \Gamma$, where $\Gamma$ is a finite subgroup of automorphisms of $T^3$. Using the classification of Seifert fiber spaces \cite{Basic3Manifold, Scott_1983}, one finds the $6$ Torus or Klein bottle bundles over the circle up to affine equivalence.
    \end{itemize}
\end{itemize}

\subsection{Spherical manifolds}
\begin{Definition}[Spherical structure]
$(X,G)$-structures are called spherical, when
$$ X = S^n = \left\{\, (x_0, x_1, \ldots, x_n)^T \,\middle|\, \sum_{i=0}^n x_i^2 = 1, x_i \in \mathbb{R} \,\right\}$$
$$ G = O(n+1),$$
where $G$ acts on $X$ by
$$ g.x = A \cdot x$$
for $g = A \in G$ and $x \in X$.
\end{Definition}

Remark. $O(n)$ preserves a Riemannian metric with sectional curvature $1$ on $S^n$ pulled back from the metric
$$ds^2 = dx_0^2 + dx_1^2 + \ldots + dx_{n-1}^2 + dx_n^2$$
on the ambient space $\mathbb{R}^{n+1}$.

What makes spherical manifolds easy to classify is that $O(n)$ is compact and hence any discrete subgroup of $O(n)$ is finite.

\subsubsection{Closed orientable examples in low dimensions}
\begin{itemize}
\item Dim = 1:
    \begin{itemize}
    \item $S^1$, the fundamental example.
    \item $S^1 = \widehat{S^1}/\Gamma$, where $\Gamma$ is a cyclic subgroup generated by a rotation.
    \end{itemize}
\item Dim = 2:
    \begin{itemize}
    \item $S^2$, the fundamental example.
    \end{itemize}
\item Dim = 3: Using the classification of discrete subgroups of $O(3)$, they are as follows.
    \begin{itemize}
    \item $S^3$, the fundamental example.
    \item $S^3/2I$, the Poincar\'{e} homology sphere ( also known as Poincar\'{e} dodecahedral space ), where $2I$ is the binary icosahedral group.
    \item $L(p; q) = S^3/\Gamma$, the lens spaces, where $\Gamma$ is isomorphic to $\mathbb{Z}/p$ and its action on $S^3$ is generated by $(z_1, z_2) \mapsto (e^{2\pi i/p}\cdot z_1, e^{2\pi i q/p}\cdot z_2)$ and $$S^3 = \left\{\, (z_1, z_2) \in \mathbb{C}^2 \,\middle|\, |z_1|^2 + |z_2|^2 = 1 \,\right\}.$$
    \end{itemize}
\end{itemize}

\subsection{Hyperbolic manifolds}
\begin{Definition}[Hyperbolic structure]
$(X,G)$-structures are called hyperbolic, when
$$ X = \mathbb{H}^n = \left\{\, (x_0, x_1, \ldots, x_n)^T \,\middle|\, \sum_{i=0}^{n-1} x_i^2 - x_n^2 = -1, x_i \in \mathbb{R}, x_n>0 \,\right\}$$
$$ G = PO(n, 1),$$
where $G$ acts on $X$ by
$$ g.x = \left\{ \begin{array}{lll} A \cdot x &,& \mbox{if } (A \cdot x)_n > 0 \\-A \cdot x &,& \mbox{if } (A \cdot x)_n < 0 \end{array} \right.$$
for $g$ represented by $A \in O(n,1)$ and $x \in X$.
\end{Definition}

Remark. $PO(n, 1)$ preserves a Riemannian metric with sectional curvature $-1$ on $\mathbb{H}^n$ pulled back from the metric
$$ds^2 = dx_0^2 + dx_1^2 + \ldots + dx_{n-1}^2 - dx_n^2$$
on the ambient space $\mathbb{R}^{n+1}$.

There are infinitely many hyperbolic surfaces in dimension 2, which demonstrates the richness of this geometry. According to Thurston's geometrization, there are enormous numbers of examples in dimension 3, yet their classification is not completely understood. Research on hyperbolic 3-manifolds remains the main stream area in 3 topology.

\subsubsection{Closed orientable examples in low dimensions}
\begin{itemize}
\item Dim = 1:
    \begin{itemize}
    \item $H^1/\Gamma$, where $\Gamma$ is a cyclic subgroup generated by a hyperbolic translation.
    \end{itemize}

\item Dim = 2:
    \begin{itemize}
    \item $\Sigma_g$, surfaces of genus $g>0$. A typical hyperbolic structure on $\Sigma_g$ can be obtained by side pairing of a regular $4g$-gon with interior angels equal to $\pi/2g$.
    \end{itemize}
    \begin{figure}[h]
    \begin{center}
    \scalebox{0.5} {\epsffile{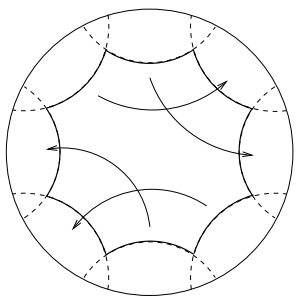}}
    \end{center}
    \caption{A hyperbolic structure on genus 2 surface in the Poincar\'{e} disk model.} \label{fig:hyperbolic_suface}
    \end{figure}
\item Dim = 3:

    We list the following examples without further defining and discussing them here. Interested readers are encouraged to follow \cite{3Manifolds, 3GeometryTopology} for details.
    \begin{itemize}
    \item Dehn surgeries on links.
    \item Most of Haken manifolds.
    \item Seifert-Weber Space.
    \end{itemize}
\end{itemize}

%% file: Chapter3.tex

\section{Affine mainfolds}\label{Ch3}
\subsection{Affine structures}
\begin{Definition}[Affine structure]
$(X,G)$-structures are called affine, when
$$ X = \mathbb{A}^n = \left\{\, (x_0, \ldots, x_{n-1})^T \,\middle|\, x_i \in \mathbb{R} \,\right\}$$
$$ G = \mbox{Aff}(n) = \mathbb{R}^n \rtimes GL(n),$$
where $G$ acts on $X$ by
$$ g.x = A \cdot x + b$$
for $g = (b, A) \in G$ and $x \in X$.
\end{Definition}

Remark. Sometimes we use the notation $(A | b)$ to denote an element $(b, A)$ in $\mbox{Aff}(n)$. Since $\mathbb{E}^n = \mathbb{A}^n$ and $E(n)$ is a subgroup of $\mbox{Aff}(n)$, $(\mathbb{E}^n, E(n))$ is a sub-geometry of $(\mathbb{A}^n, \mbox{Aff}(n))$. Hence all Euclidean manifolds are affine.

As mentioned before in \S\ref{subsec:examples}, affine structures need not be complete. There are 2 major types and also other examples:
\begin{itemize}
\item Complete Case - The space forms: $M = \mathbb{A}^n/\Gamma$.

    Here $\Gamma \subset \mbox{Aff}(n)$ is a discrete subgroup that acts freely and properly discontinuously on $X$.

    Two famous conjectures about complete affine structures are:
    \begin{Conjecture}[Markus conjecture]
    An affine structure on a closed manifold is complete if and only if the holonomy group preserves a constant volume form, i.e. it is contained in $\mathbb{R}^n \rtimes SL(n)$.
    \end{Conjecture}

    \begin{Conjecture}[Auslander conjecture]
    The fundamental group of a closed manifold which admits a complete affine structure is virtually polycyclic.
    \end{Conjecture}

\item Radiant Case: $(\mathbb{A}^n, GL(n))$-structures.
    \begin{Definition}[Radiant structures]
    An affine structure is called radiant, if the holonomy group fixes a point in $\mathbb{A}^n$, i.e. it is conjugate to an $(\mathbb{A}^n, GL(n))$-structure.
    \end{Definition}

    An important geometric invariant of radiant structures is the radiant vector field.
    \begin{Definition}[Radiant vector field]
    The vector field $\sum\limits_{i=0}^{n-1} x^i \frac{\partial}{\partial x^i}$ on $\mathbb{A}^n$ is invariant under $GL(n)$, and hence it descends to a vector field on any radiant manifolds.
    \end{Definition}

    A basic example of radiant manifold is a Hopf manifold: $$M \, = \, \mathbb{A}^n-\{0\} / \langle x \sim 2x \rangle \, \simeq \, S^{n-1}\times S^1.$$

\item Other Cases:
    There exist other affine structures that are neither complete nor radiant.
\end{itemize}

\subsection{Closed orientable affine manifold in dimension 2}
Benz\'{e}cri showed the following result \cite{Benzecri_1960}. It also follows from Milnor's more algebraic approach on flat bundles \cite{Milnor_1958}.
\begin{Theorem}
Let $M$ be a closed 2-dimensional affine manifold. Then $\chi(M) = 0$.
\end{Theorem}

Nagano and Yagi classified the affine structures on the real two-torus \cite{Nagano-Yagi_1973} ( See also\cite{Baues-Goldman_2005, Benoist_2000} ).

The holonomy of an affine structure on $T^2$ is a representation $$\pi_1(M) = \mathbb{Z}^2 \to \mbox{Aff}(2),$$ which can be extended to $$\mathbb{R}^2 \to \mbox{Aff}(2)$$ by the property of nilpotent groups \cite{Discrete}. One can show that any homomorphism from $\mathbb{R}^2$ to $\mbox{Aff}(2)$ with a 2-dimensional image is conjugate in $\mbox{Aff}(2)$ to one of the six in the following. Then we can recover an affine structure on the torus by choosing a lattice in $\mathbb{R}^2$.

Remark. The figures displayed below are just one of the examples by choosing a particular lattice $L \doteq \mathbb{Z}^2$. There are many interesting pictures of tessellation from different choices of lattices ( See \cite{Baues-Goldman_2005, Benoist_2000} ).

Recall that we use the notation $\left( A \middle| b \right)$ to represent an affine transformation with its linear part equal to $A$ and translational part equal to $b$.

\begin{itemize}
\newpage
\item $\phi_1: (s,t) \mapsto \left(\begin{array}{cc} 1 & 0 \\ 0 & 1\end{array} \middle| \begin{array}{c} s \\ t\end{array} \right)$ where $(s, t) \in \mathbb{R}^2$. A lattice of $\mathbb{R}^2$ can act on $\mathbb{A}^2$ freely and properly discontinuously via $\phi_1$. This structure is complete. It is the one conjugate to a Euclidean structure.

    \begin{figure}[h]
    \begin{center}
    \scalebox{0.5} {\epsffile{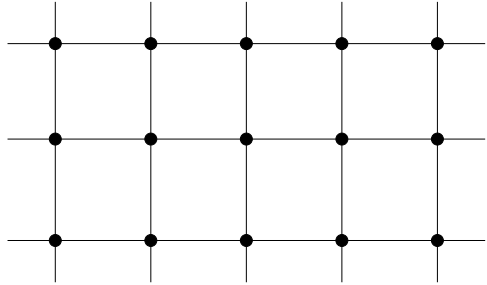}}
    \end{center}
    \caption{An affine structure from a lattice for $\phi_1$.} \label{fig:T2_1}
    \end{figure}

\bigskip
\bigskip

\item $\phi_2: (s,t) \mapsto \left(\begin{array}{cc} 1 & t \\ 0 & 1\end{array} \middle| \begin{array}{c} s + \frac{t^2}{2} \\ t\end{array} \right)$ where $(s,t) \in \mathbb{R}^2$. A lattice of $\mathbb{R}^2$ can act on $\mathbb{A}^2$ freely and properly discontinuously via $\phi_2$. This structure is complete.

    \begin{figure}[h]
    \begin{center}
    \scalebox{0.5} {\epsffile{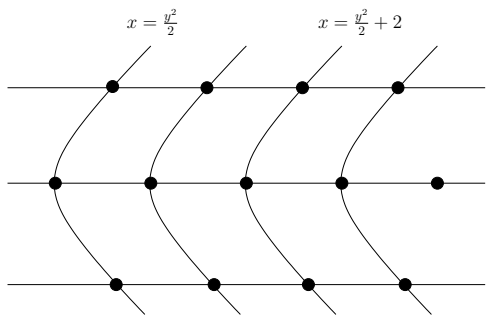}}
    \end{center}
    \caption{An affine structure from a lattice for $\phi_2$.} \label{fig:T2_2}
    \end{figure}

    Remark. In the figure above, we can identify $\mathbb{R}^2$ with $\mathbb{A}^2$ by letting $\phi_2(s,t)$ act on $(0, 0)^T$, which is given by $(s,t) \mapsto (x, y) = (s + t^2/2, t)$.

\newpage
\item $\phi_3: (s,t) \mapsto \left(\begin{array}{cc} \exp(s) & \exp(s)t \\ 0 & \exp(s) \end{array} \middle| \begin{array}{c} 0 \\ 0\end{array} \right)$ where $(s,t) \in \mathbb{R}^2$. A lattice of $\mathbb{R}^2$ can act on the open upper half plane freely and properly discontinuously via $\phi_3$. This structure is radiant.

    \begin{figure}[h]
    \begin{center}
    \scalebox{0.5} {\epsffile{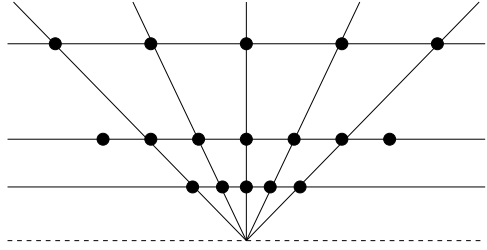}}
    \end{center}
    \caption{An affine structure from a lattice for $\phi_3$.} \label{fig:T2_3}
    \end{figure}

\bigskip
\bigskip

\item $\phi_4: (s,t) \mapsto \left(\begin{array}{cc} \exp(s) & 0 \\ 0 & \exp(t) \end{array} \middle| \begin{array}{c} 0 \\ 0\end{array} \right)$ where $(s,t) \in \mathbb{R}^2$. A lattice of $\mathbb{R}^2$ can act on the open first quadrant freely and properly discontinuously via $\phi_4$. This structure is radiant.

    \begin{figure}[h]
    \begin{center}
    \scalebox{0.5} {\epsffile{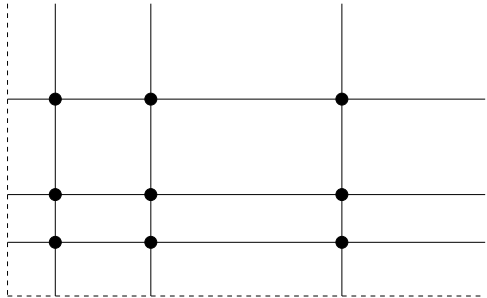}}
    \end{center}
    \caption{An affine structure from a lattice for $\phi_4$.} \label{fig:T2_4}
    \end{figure}

\newpage
\item $\phi_5: (s,t) \mapsto \left(\begin{array}{cc} \exp(s)\cos(t) & \exp(s) \sin(t) \\ -\exp(s) \sin(t) & \exp(s) \cos(t) \end{array} \middle| \begin{array}{c} 0 \\ 0\end{array} \right)$ where $(s,t) \in \mathbb{R}^2$. A lattice of $\mathbb{R}^2$ can act on the universal cover of $\mathbb{A}^n - \{0\}$ freely and properly discontinuously via $\phi_5$. This structure is radiant.

    \begin{figure}[h]
    \begin{center}
    \scalebox{0.5} {\epsffile{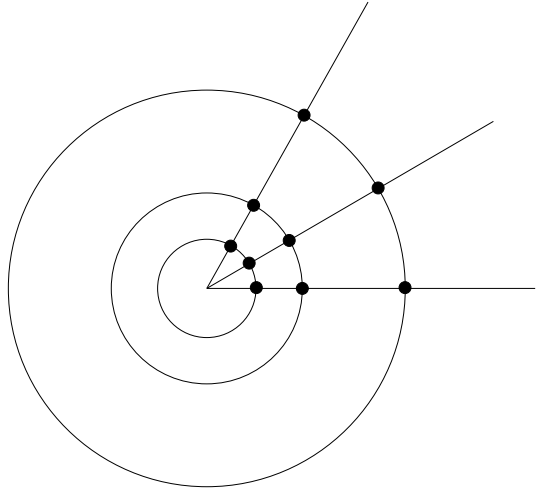}}
    \end{center}
    \caption{An affine structure from a lattice for $\phi_5$.} \label{fig:T2_5}
    \end{figure}

\bigskip
\bigskip

\item $\phi_6: (s,t) \mapsto \left(\begin{array}{cc} 1 & 0 \\ 0 & \exp(t) \end{array} \middle| \begin{array}{c} s \\ 0\end{array} \right)$ where $(s,t) \in \mathbb{R}^2$. A lattice of $\mathbb{R}^2$ can act on the open upper half plane freely and properly discontinuously via $\phi_6$. This structure is neither complete nor radiant.

    \begin{figure}[h]
    \begin{center}
    \scalebox{0.5} {\epsffile{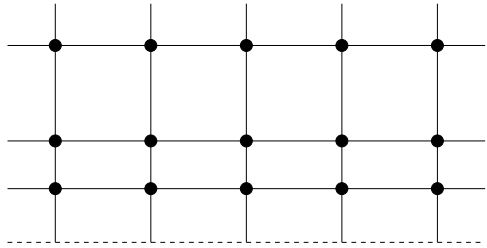}}
    \end{center}
    \caption{An affine structure from a lattice for $\phi_6$.} \label{fig:T2_6}
    \end{figure}
\end{itemize}

\subsection{Closed orientable affine manifolds in dimension 3}
Unlike the case of dimension 2, we do not have a list of all closed orientable affine 3-manifolds. Contrary to the Bieberbach Theorem ( Theorem \ref{thm:Bieberbach} ) of Euclidean manifolds, Auslander found that there exists a family of countably many distinct affine 3-manifolds. Nevertheless the affine manifolds with complete structures have been classified \cite{Fried-Goldman_1983}, so have those with radiant structures \cite{Choi_decomp_2001}. We also have affine structures on a 3-manifold, which is neither complete nor radiant.

\begin{Definition}[Prime 3-manifolds]
A 3-manifold is prime if it cannot be expressed as a non-trivial connected sum of two 3-manifolds.
\end{Definition}

We note that all known affine 3-manifolds below are prime. This fact motivates us to ask the following question, which is the special case in dimension 3 of the basic topological question from the introduction.
\begin{center}\it ``Are all closed orientable affine 3-manifolds prime?''\end{center}

Since we have sphere theorem and other techniques ready from basic 3-manifold topology, this problem might be attackable. We do not have the answer yet, however let us state it as a conjecture.
\begin{Conjecture}
All closed orientable affine 3-manifolds are prime.
\end{Conjecture}

Let us take a look at what affine 3-manifolds we have known in the following.

\subsubsection{Complete structures}

Fried and Goldman classified complete closed affine manifolds as follow \cite{Fried-Goldman_1983}.\label{thm:complete}
\begin{Theorem}
Let $M$ be a closed 3-manifold. Then the following are equivalent:
\begin{itemize}
  \item $M$ admits a complete affine structure.
  \item $\pi_1(M)$ is solvable and $M$ is aspherical.
  \item $M$ is finitely covered by a 2-torus bundle over the circle.
\end{itemize}
\end{Theorem}

The first step of their proof is to show that the holonomy group is virtually solvable, i.e. it contains a solvable subgroup of finite index.

First of all, an important observation of a complete structure is that the linear part of any element in the holonomy group must have $1$ as an eigenvalue, for otherwise it has a fixed point and can not act freely.

Therefore we have
\begin{Lemma}[Lemma 2.3 in \cite{Fried-Goldman_1983}]
The holonomy group of a complete affine structure is contained in the algebraic subgroup of $\mbox{Aff}(n)$ consisting of those elements with the linear part containing $1$ as an eigenvalue.
\end{Lemma}

Second, we note
\begin{Lemma}[Lemma 2.6 in \cite{Fried-Goldman_1983}]
The only semisimple connected subgroups of $SL(3, \mathbb{R})$ are
$$I \mbox{, } SO(3) \mbox{, } SO(2,1)^0 \mbox{, } SL(2,\mathbb{R})\times\{1\} \mbox{ and } SL(3, \mathbb{R}).$$
\end{Lemma}

Using the above lemmas and the fact that the holonomy group acts cocompactly, one can rule out all the non-trivial subgroups listed in the above lemma as the Levi factor of the subgroup of the holonomy group, which consists of those elements in $SL(3, \mathbb{R})$. ( A Levi factor is the semisimple part of the decomposition of a group into a semidirect product of a semisimple subgroup and the maximal solvable subgroup. )

Hence we have
\begin{Theorem}[Theorem 2.1 in \cite{Fried-Goldman_1983}]
The holonomy group of a complete affine structure on a closed 3-manifold is virtually solvable.
\end{Theorem}

Lastly let us define crystallographic hull as follow
\begin{Definition}[Theorem on pg.5 in \cite{Fried-Goldman_1983}]
Let $\Gamma \subset \mbox{Aff}(n)$ be virtually solvable and suppose $\Gamma$ acts properly discontinuously on $\mathbb{A}^n$. A subgroup $H \subset \mbox{Aff}(n)$ containing $\Gamma$ is called the crystallographic hull for $\Gamma$ if $H$ satisfies the following:
\begin{itemize}
    \item $H$ has finitely many components and each component of $H$ meets $\Gamma$;
    \item $H / \Gamma$ is compact;
    \item $H$ and $\Gamma$ have the same algebraic hull in $\mbox{Aff}(n)$;
    \item Every isotropy group of $H$ on $\mathbb{A}^n$ is finite.
\end{itemize}
\end{Definition}

By further study on the crystallographic hull for the holonomy group, they show that the holonomy group is actually solvable. Hence the fundamental group of a closed complete affine manifold is solvable. The classification of crystallographic hull also gives
\begin{Theorem}[Theorem on pg.21 in \cite{Fried-Goldman_1983}]
A closed complete affine 3-manifold is finitely covered by a Solvmanifold with a left invariant affine structure.
\end{Theorem}

On one hand, any complete affine 3-manifold is covered by $\mathbb{A}^3$ by definition, and hence is aspherical. Then it is finitely covered by a 2-torus bundle over the circle by a theorem of Evans-Moser \cite{Evans_Moser_1972}.

On the other hand, suppose $M$ is finitely covered by $N$ which is a 2-torus bundle over the circle. If we identify $T^2$ with $\mathbb{R}^2 / \mathbb{Z}^2$, a homeomorphism of $T^2$ is homotopic to one induced from an $SL(2,\mathbb{Z})$ action on $\mathbb{R}^2$, and hence $N$ admits an affine structure. Any deck transformation of the covering $N \to M$ can also be homotopic to one induced from an affine action on $\mathbb{R}^2$. Therefore $M$ is homotopy equivalent to a complete affine manifold. Then one can show that a homotopy-equivalence is homotopic to a homeomorphism in this case ( see \cite{Fried-Goldman_1983} for more details ).

\subsubsection{Radiant structures}
Barbot and Choi classified radiant closed affine manifolds as follow \cite{Choi_decomp_2001}.
\begin{Theorem}
A closed 3-manifold admits a radiant affine structure if and only if it is one of the following:
\begin{itemize}
\item Benz\'{e}cri suspension over $\Sigma_g$ with $g>1$.
\item Generalized affine suspension of $\Sigma_g$ with $g \leq 1$.
\item Finitely covered by a 2-torus bundle over the circle.
\end{itemize}
\end{Theorem}

We postpone our discussion of affine suspension later. The ``if'' direction is straightforward, while the ``only if'' direction follows from the fact that any compact radiant affine 3-manifold admits a total cross-section to the radiant flow generated by the radiant vector field, which itself follows from a more interesting geometric decomposition as follow.
\begin{Theorem} \label{thm:radiantdecomp}
Let $M$ be a compact radiant affine 3-manifold with empty or totally geodesic boundary. Then $M$ decomposes along the union of finitely many disjoint totally geodesic tori or Klein bottles, tangent to the radial flow, into a disjoint union of
\begin{itemize}
  \item Concave affine 3-manifolds.
  \item Crescent cone affine 3-manifolds.
  \item Convex affine 3-manifolds.
\end{itemize}
\end{Theorem}

We need the following definitions before proceeding to talk about the decomposition.

\begin{figure}[h]
\begin{center}
\scalebox{0.6} {\epsffile{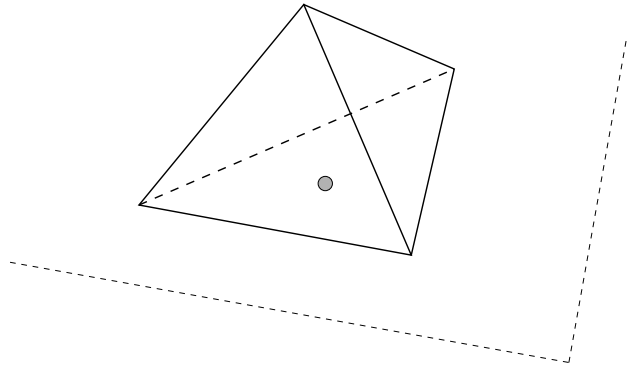}}
\end{center}
\caption{The closed upper half space of $\mathbb{A}^3$ with a point removed is not 2-convex.} \label{fig:not_2_convex}
\end{figure}

\begin{Definition}[$m$-convex]
$M$ is said to be $m$-convex, if it satisfies the following property: for any subset $S$ of $\hat{M}$, if $S$ is homeomorphic via $dev$ to a $m$-simplex without the interior of only one face, then there exists a subset $S'$ of $\hat{M}$ containing $S$ such that $S'$ is actually homeomorphic via $dev$ to the $m$-simplex.
\end{Definition}

Remark. 1-convexity is equivalent to the usual convexity. Note also that m-convexity implies (m+1)-convexity.

\begin{Definition}[Crescent]
A crescent is a closed subset $C$ of $\hat{M}$ such that $dev$ is a homeomorphism from $C$ onto $dev(C)$, such that the interior of $dev(C)$ is an open half space of $\mathbb{A}^3$.
\end{Definition}

\begin{Definition}[Equivalence of crescent]
Two crescents $C$ and $C'$ are said to be equivalent, if there exists a chain of crescents $C_1=C, C_2, \ldots, C_k=C'$ such that $C_i \cap C_{i+1} \neq \emptyset$ for $i = 1, \ldots, k-1$.
\end{Definition}

\begin{Definition}[Concave affine 3-manifolds]
$M$ is said to be a concave affine 3-manifold, if $\hat{M}$ is a union of crescents.

\end{Definition}
\begin{figure}[h]
\begin{center}
\scalebox{0.6} {\epsffile{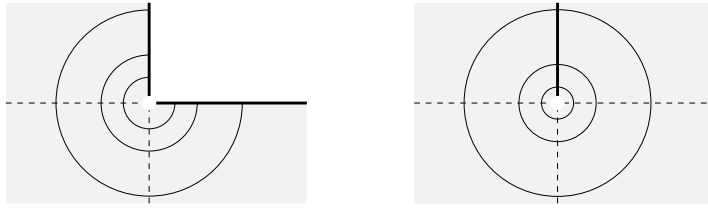}}
\end{center}
\caption{Examples of concave affine manifolds as unions of crescents.} \label{fig:concave}
\end{figure}

Remark, we show two examples of concave affine manifolds in Figure \ref{fig:concave}. The shaded region of the left one is a union of two crescents ( the left half space and the bottom half space ), while the shaded region of the right one is a union of three crescents ( the left half space, the bottom haft space and the right half space ).

\begin{Definition}[Crescent cone 3-manifolds]
$M$ is said to be a crescent cone 3-manifolds, if $dev$ is a homeomorphism from $\hat{M}$ onto $dev(\hat{M})$ such that the interior of $dev(\hat{M})$ is the intersection of two open half spaces of $\mathbb{A}^3$.
\end{Definition}

\begin{figure}[h]
\begin{center}
\scalebox{0.6} {\epsffile{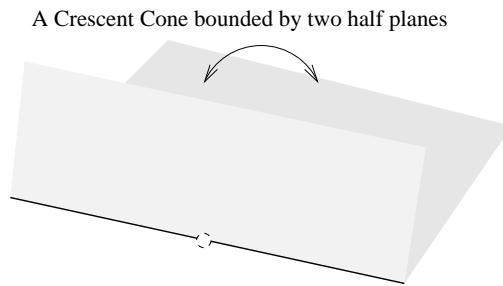}}
\end{center}
\caption{An example of a crescent cone, which is 2-convex but not convex.} \label{fig:crescent_cone}
\end{figure}

The decomposition in Theorem \ref{thm:radiantdecomp} is done in the following way:
\begin{itemize}
\item If $M$ is not 2-convex, Choi showed that there exist crescents ( Figure \ref{fig:not_2_convex} ). Since the union of the crescents in an equivalence class of crescents is a concave affine submanifold, we attain the first decomposition of $M$ into:
    \begin{itemize}
      \item Concave affine 3-manifolds.
      \item $2$-convex affine 3-manifolds.
    \end{itemize}

\item Then if a 2-convex affine 3-manifold is not convex, Choi showed that there exist crescent cones ( Figure \ref{fig:crescent_cone} ). Hence we can decompose all 2-convex affine 3-manifolds into:
    \begin{itemize}
      \item Crescent cone affine 3-manifolds.
      \item Convex affine 3-manifolds.
    \end{itemize}
\end{itemize}
Now the decomposition has been completed.

\subsubsection{Projective structures}
Affine suspension is a construction of a radiant affine manifold from a projective manifold of one dimension lower. We first need to discuss projective structures.

\begin{Definition}[Projective structure]
$(X,G)$-structures are called projective, when
$$ X = \mathbb{RP}^n = \mathbb{R}^n-\{0\} / \langle x \sim \lambda x, \lambda \neq 0\rangle$$
$$ G = PGL(n+1),$$
where $G$ acts on $X$ by
$$ g.[x] = [A \cdot x]$$
for $g$ represented by $A \in GL(n+1)$ and $x \in \mathbb{R}^n-\{0\}$.
\end{Definition}

Remark. As mention before, Euclidean, affine, spherical and hyperbolic geometry are all sub-geometries of projective geometry. Their relations can be seen in Figure \ref{fig:xg_hier}.

The study of projective geometry is largely motivated by hyperbolic geometry. Since the developing image of a hyperbolic surface is inside the light cone, which is a convex disk when viewing projectively. This motivates us to consider a sub-category of projective structures, which is called convex projective structures \cite{Goldman_1990}.

\begin{Definition}
A projective structure on $M$ is called convex, if $dev$ is a homeomorphism from $\hat{M}$ onto a convex domain of $\mathbb{RP}^n$.
\end{Definition}

We also like to mention the definitions of marked structures and deformation space here.

\begin{Definition}[Marked structures]
A marked $\mathbb{RP}^n$-structure on $M$ is given by a pair $(\phi, N)$ such that $N$ is a $\mathbb{RP}^n$-manifold and $\phi$ is a diffeomorphism from $M$ to $N$.
\end{Definition}

\begin{Definition}[Deformation space]
The deformation space of marked $\mathbb{RP}^n$-structures on $M$ is the quotient space of all the marked $\mathbb{RP}^n$-structures on $M$ modulo the equivalence relation defined by the following:
$$ (\phi, N) \sim (\phi', N')$$
if $\phi' \circ \phi^{-1}: N \to N'$ is isotopic to a projective map from $N$ to $N'$.
\end{Definition}

The main result in \cite{Goldman_1990} is
\begin{Theorem}
The deformation space of marked convex projective structures on a surface $\Sigma_g$ for $g>1$ is diffeomorphic to  an open cell of dimension $16(g-1)$.
\end{Theorem}

Combining Choi's result on the decomposition of a projective manifold along some totally geodesic boundary,
Choi and Goldman classified all the real projective structures on closed surfaces \cite{Choi_Goldman_1997}.
\begin{Theorem}
The deformation space of marked projective structures on a surface $\Sigma_g$ for $g>1$ is diffeomorphic to a countable union of open cells of dimension $16(g-1)$.
\end{Theorem}
What we need from the theorem is that any projective structure on $\Sigma_g$ for $g>1$ can be constructed from a convex projective structure by grafting projective annuli. Then by affine suspension, we can construct some exotic radiant affine structures on $\Sigma_g \times S^1$.

\subsubsection{Example of grafting a projective surface}
The following example illustrate this process.
\begin{itemize}
\item Let us take a closed hyperbolic surface $\Sigma_g$. Then its developing image is the disk inside the light cone projectively. See Figure \ref{fig:axes}, the darker shaded region is the developing image.

\item Take any close geodesic $\gamma$ on $\Sigma_g$, and cut $\Sigma_g$ open along $\gamma$ into $\Sigma_g \backslash \gamma$. Then $\gamma$ can represent a deck transformation of $\Sigma_g$. We use the same notation for the corresponding element in the fundamental group.

\item Then $hol(\gamma)$ is a hyperbolic element in $SO(2,1) \cong PO(2,1)$, hence it has three axes $v^+$, $v^-$ and $v^0$ corresponding to eigenvectors of eigenvalues $\lambda$, $\lambda^{-1}$ and $1$ respectively, where $\lambda > 1$. See Figure \ref{fig:axes} and \ref{fig:axes_plane}.
\begin{figure}[h]
\begin{center}
\scalebox{0.6} {\epsffile{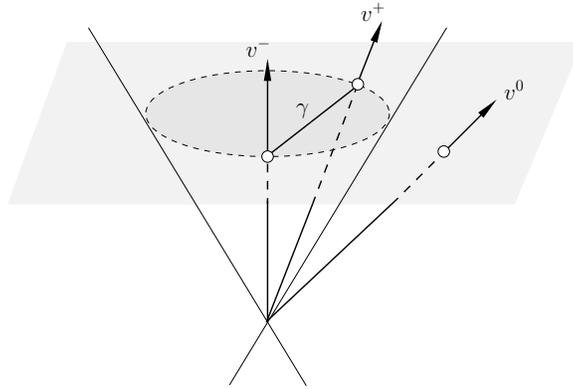}}
\end{center}
\caption{Three axes of $hol(\gamma)$ and a plane.} \label{fig:axes}
\end{figure}
\begin{figure}[h]
\begin{center}
\scalebox{.565} {\epsffile{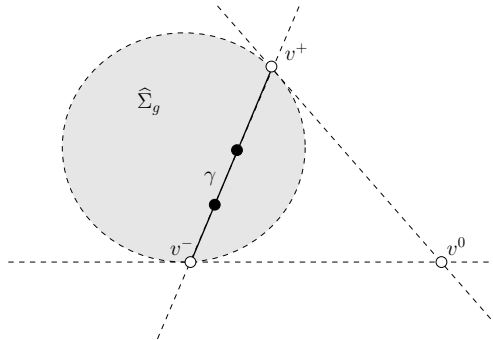}}
\end{center}
\caption{Three axes of $hol(\gamma)$ seen on a plane.} \label{fig:axes_plane}
\end{figure}
\item One can find a projective structure ( many ) on the annulus $A \doteq \gamma \times [0,1]$ with the holonomy group generated exactly by $hol(\gamma)$.

\begin{figure}[h]
\begin{center}
\scalebox{0.565} {\epsffile{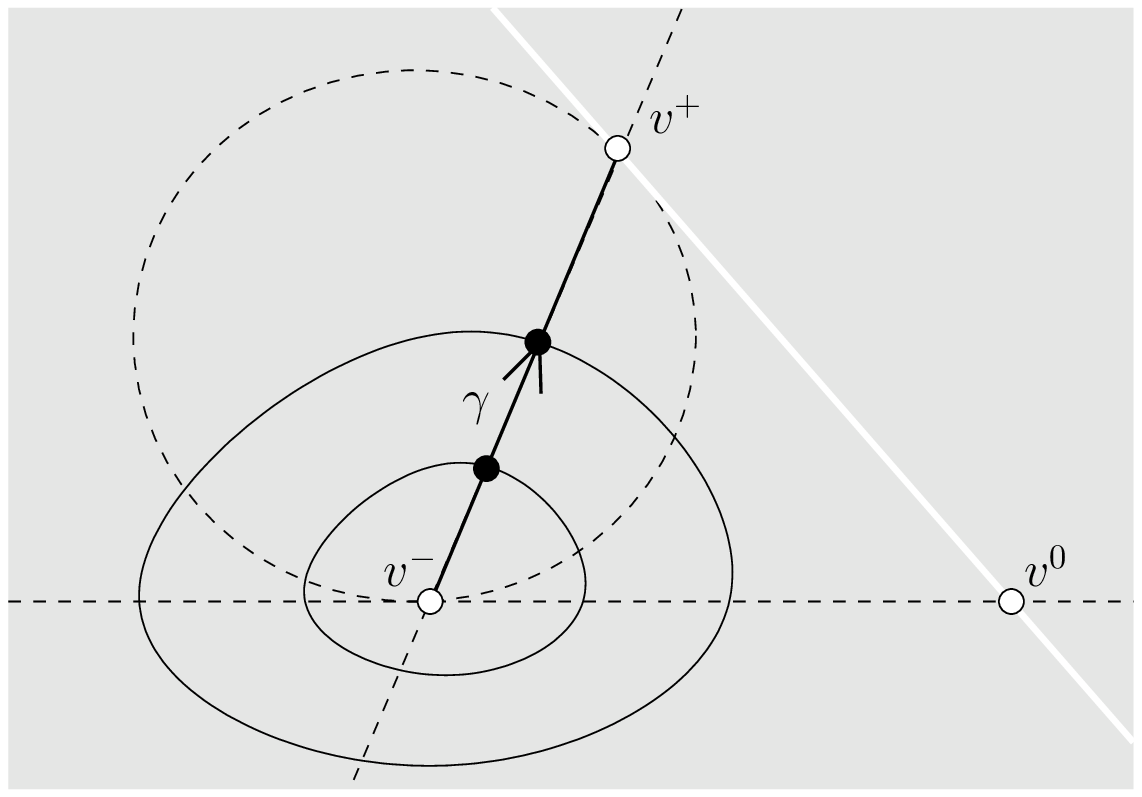}}
\end{center}
\begin{center}
\scalebox{0.565} {\epsffile{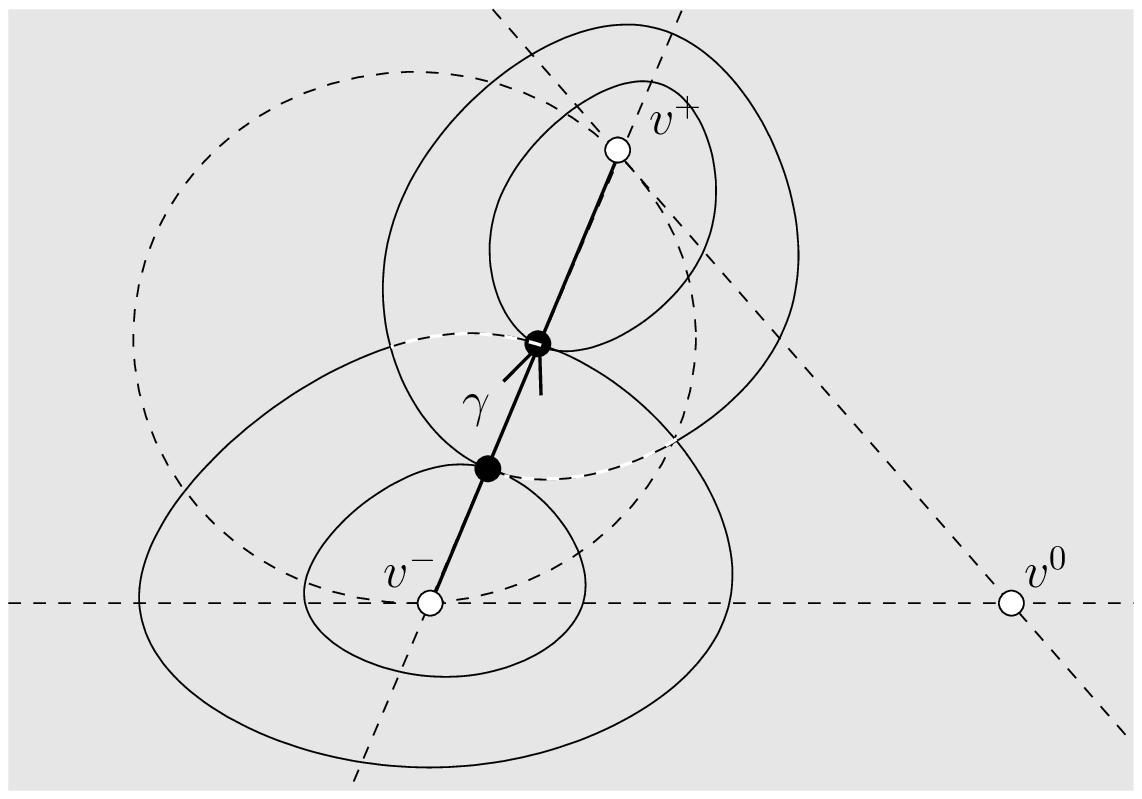}}
\end{center}
\caption{Two projective annuli with holonomy generated by $hol(\gamma)$.} \label{fig:proj_annuli}
\end{figure}

See Figure \ref{fig:proj_annuli} for two examples of projective structures on the annulus $A$. The shaded regions are the developing images: the developing image of the top one is $\mathbb{RP}^2$ minus the point $v^-$ and the line through $v^+$ and $v^0$; the developing image of the bottom one is $\mathbb{RP}^2$ minus the three points $v^-$, $v^+$ and $v^0$.

\item Now the projective structures on $\Sigma_g \backslash \gamma$ and $A$ are exactly the same along the ends ( the developing maps restricted to a small neighborhood of their boundary are the same ), hence we can concatenate the ends ( so called grafting, see Figure \ref{fig:proj_concatenation} ) and get a new projective surface, which is homeomorphic to $\Sigma_g$.
\begin{figure}[h]
\begin{center}
\scalebox{0.6} {\epsffile{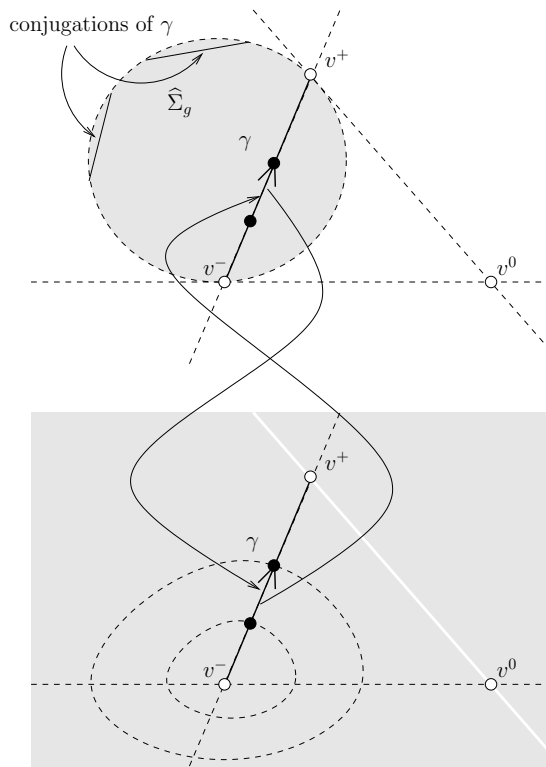}}
\end{center}
\caption{A new projective structure from grafting.} \label{fig:proj_concatenation}
\end{figure}
\end{itemize}
Remark. The projective structure constructed above from grafting has the property that the developing map is surjective onto $\mathbb{RP}^2$ but it is not a covering map.

\subsubsection{Affine suspension}
Now given a projective surface $\Sigma$, we can choose a developing pair $dev$ and $hol$. Note that $S^2$ is the universal cover of $\mathbb{RP}^2$, $dev$ actually factors through $S^2$. Let us still use $dev$ to denote this map $$ dev: \hat{\Sigma} \to S^2.$$ We can now pull back the canonical trivial line bundle over $S^2$, which is just $\mathbb{R}^3-\{0\}$, to $\hat{\Sigma}$. Hence we have a map from the trivial bundle to $\mathbb{R}^3-\{0\}$:
$$ \widetilde{dev}: \hat{\Sigma} \times \mathbb{R} \to \mathbb{R}^3-\{0\}$$
$$ (\hat{m}, t) \mapsto \exp(t) dev(\hat{m}),$$
and $$hol: \pi_1(\Sigma) \to PGL(3) \cong SL^{\pm}(3) \subset GL(3).$$
Let us use $\widetilde{hol}$ to denote this homomorphism $$\widetilde{hol}: \pi_1(\Sigma) \to SL^{\pm}(3),$$
which makes $\pi_1(\Sigma)$ act on $\mathbb{R}^3-\{0\}$.

Therefore there is a natural $\pi_1(\Sigma)$ action on the pull back bundle $\hat{\Sigma} \times \mathbb{R}$, which still acts freely and properly discontinuously and makes $(\widetilde{dev}, \widetilde{hol})$ a developing pair for a $(\mathbb{R}^3-\{0\}, SL^{\pm}(3))$-structure, a radiant affine structure.

In order to obtain a compact manifold, note that $\pi_1(\Sigma)$ acts properly discontinuously of $\hat{\Sigma}$, hence if we choose $\lambda > 0$ large enough and let $\mathbb{Z}$ be the cyclic group generated by this homothety $\lambda\cdot I$. Then $\pi_1(\Sigma) \times \mathbb{Z}$ acts on $\hat{\Sigma} \times \mathbb{R}$ freely, properly discontinuously and cocompactly.

\begin{Definition}[Benz\'{e}cri suspension]
The process described above to construct a radiant affine manifold from a projective manifold of one dimension lower together with a homothety action on the radial direction is called a Benz\'{e}cri suspension.
\end{Definition}

More generally, instead of using a homothety, if we have a projective automorphism $\phi$ of $\Sigma$, i.e. there exists $\psi \in PGL(3)$ and a homomorphism $$ \Phi: \pi_1(\Sigma) \to \pi_1(\Sigma),$$such that
$$ dev(\phi.\hat{x}) = \psi.dev(\hat{x})$$
$$ \phi.\gamma.\hat{x} = \Phi(\gamma).\phi.\hat{x}$$
for all $\hat{x} \in \hat{\Sigma}$ and $\gamma \in \pi_1(\Sigma)$. We can lift $\psi$ to $GL(3)$ with the absolute value of the determinant large enough, and then the semidirect product $$\pi_1(\Sigma) \rtimes_{\Phi} \langle \phi \rangle$$ acts on $\hat{\Sigma} \times \mathbb{R}$ freely, properly discontinuously and cocompactly.

\begin{Definition}[Affine suspension]
The process described above to construct a radiant affine manifold from a projective manifold of one dimension lower together with a projective automorphism of the underlying projective manifold is called an affine suspension.
\end{Definition}

\subsubsection{Other affine structures}
The developing images of the non-complete examples of affine structures ( non-complete affine structures on $T^2$, radiant structures on affine 3-manifolds ) we have seen so far are cones. Goldman gave two examples to show that this is not the case in general\cite{Goldman_1981}, and in particular they are not radiant either.

Let $$0_+ = \left\{\, (x,y,z)\in \mathbb{R}^3 \,\middle|\, y^2-2x > 0 \,\right\}$$
$$0_- = \left\{\, (x,y,z)\in \mathbb{R}^3 \,\middle|\, y^2-2x < 0 \,\right\},$$
and let $G$ be the group of the form
$$ \left( \begin{array}{ccc} e^{2s} & e^s t & 0\\
0 & e^s & 0 \\ 0 & 0 & e^{-s}
\end{array}\middle| \begin{array}{c}
\frac{t^2}{2}\\t\\u
\end{array}\right),$$
where $s, t, u \in \mathbb{R}$.

Then $G$ acts simply transitively on both $0_+$ and $0_-$. Hence if we choose a discrete cocompact subgroup $\Gamma$ of $G$, then $0_+/\Gamma$ and $0_-/\Gamma$ are closed affine 3-manifolds. Obviously neither $0_+$ nor $0_-$ is a cone, and $0_+$ is not even convex.

%% file: Chapter4.tex

\section{Embedded spheres of affine manifolds}\label{Ch4}
From all the examples we have seen so far in \S3, they are either covered by $S^2 \times \mathbb{R}$ ( affine suspension of $S^2$ ), or covered by an open 3-cell ( complete affine manifolds, torus/Klein bottle bundle over a circle, affine suspension of a surface $\Sigma_g$ for $g \geq 0$, Goldman's examples ). Therefore all known closed orientable affine 3-manifolds are prime. We are really interested in the following question:
{\begin{center}`` Are all closed orientable affine 3-manifolds prime? ''\end{center}}
As mentioned before in \S3.3 we do not have the answer to the question yet. Instead, what we are going to prove in this Chapter is\\
\textbf{Theorem \ref{thm:main}}
{\it For $n \geq 3$, let $(M,\partial)$ be a compact affine $n$-manifold with boundary $\partial \simeq S^{n-1}$ and $dev: \hat{M} \to \mathbb{A}^n$ be a developing map. If
\begin{itemize}
\item $dev$ restricted to some lift $\hat{\partial}$ of $\partial$ is an embedding,
\item $dev$ maps a neighborhood of $\hat{\partial}$ to the closure of the bounded part of $\mathbb{A}^n \backslash dev(\hat{\partial})$,
\end{itemize}
then $(M, \partial)$ is homeomorphic to $(D^n, S^{n-1})$.}

Recall the following theorem from Mazur, Morse and Brown \cite{Mazur_1959, Morse_1960, Brown_1960}.
\begin{Lemma}[Generalized Schoenflies theorem]
If an imbedding $\psi: S^{n-1} \to S^n$ has the property that there is a $\phi: S^{n-1} \times [-1,1] \to S^n$ such that $\phi$ is a homeomorphism onto its image and $\phi(S^{n-1} \times \{0\}) = \psi(S^{n-1})$, then the closure of each components of $S^n - \psi(S^{n-1})$ is topologically an $n$-cell, i.e. the pair $(\psi(S^{n-1}), S^n)$ is homeomorphic to $(S^{n-1}, S^n)$.
\end{Lemma}

\begin{figure}[h]
\begin{center}
\scalebox{0.35} {\epsffile{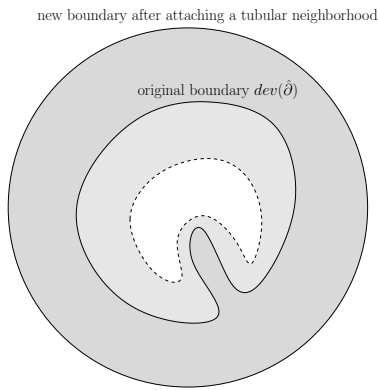}}
\end{center}
\caption{Deform the boundary to a standard sphere.} \label{fig:extend_bdry}
\end{figure}

Therefore by attaching a collar neighborhood of the boundary, we can reduce Theorem \ref{thm:main} to the following special case, in which the strictly convexity of the sphere enables us to use our geometric approach on dome bodies defined in the next subsection.
\begin{Theorem} \label{thm:convex_main}
For $n \geq 3$, let $(M,\partial)$ be a compact affine $n$-manifold with boundary $\partial \simeq S^{n-1}$ and $dev: \hat{M} \to \mathbb{A}^n$ be a developing map. If
\begin{itemize}
\item $dev$ restricted to some lift $\hat{\partial}$ of $\partial$ is an embedding onto a standard sphere which bounds a strictly convex closed solid ball in $\mathbb{A}^n$,
\item $dev$ maps a neighborhood of $\hat{\partial}$ to the closure of the bounded part of $\mathbb{A}^n \backslash dev(\hat{\partial})$,
\end{itemize}
then $(M, \partial)$ is homeomorphic to $(D^n, S^{n-1})$.
\end{Theorem}

Remark. In this case, the image under $dev$ of every component of $\partial\hat{M}$ is affinely equivalent to a standard sphere.

\subsection{Dome body}
From now on, we are working on a compact affine $n$-manifold $M$ with boundary $\partial \simeq S^{n-1}$ such that $dev(\hat{\partial})$ is affinely equivalent to a standard sphere for any component $\hat{\partial}$ of $\partial\hat{M}$.

We now fix an arbitrary translational invariant metric on $\mathbb{A}^n$. So whenever we talk about length and volume, we refer to this metric.

We will also talk about convergence of a sequence of subspaces ( e.g. lines, hyperplanes ) in the following. All the subspaces in such a sequence, of which we are going to take the limit, will always meet some common compact region. Then we just consider the Hausdorff distance between them restricted to the common compact region.

\begin{figure}[h]
\begin{center}
\scalebox{0.6} {\epsffile{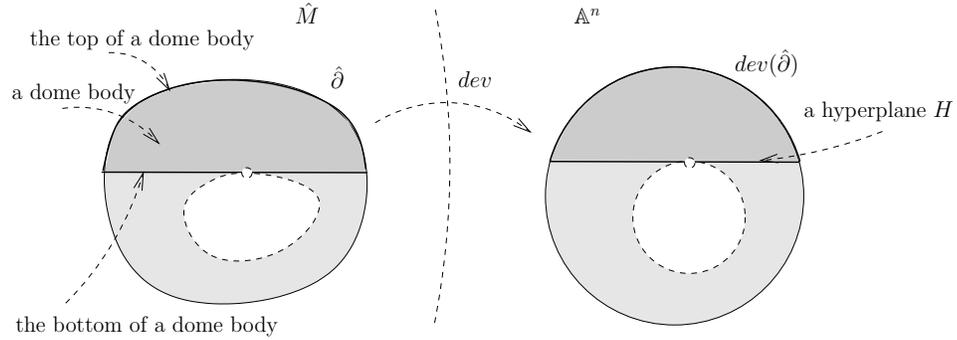}}
\end{center}
\caption{A dome body.} \label{fig:domebody}
\end{figure}
\begin{Definition}[Dome body]\label{dfn:domebody}
A subset $D$ of $\hat{M}$ is called a dome body, if \begin{itemize}
\item $D$ is the closure of its interior in $\hat{M}$,
\item $D \cap \hat{\partial} \neq \emptyset$ for some component $\hat{\partial}$ of $\partial\hat{M}$,
\item $dev$ restricted to the interior of $D$ is an embedding onto the open solid semi ball, which is the intersection of an open half space of $\mathbb{A}^n$ bounded by a hyperplane $H$ and the open ball bounded by $dev(\hat{\partial})$.
\end{itemize}
\end{Definition}

Let us define the following terms for our convenience.

\begin{Definition}[Top of a dome body]
The subset $D \cap \hat{\partial}$ is called the top of the dome body $D$.
\end{Definition}

\begin{Definition}[Bottom of a dome body]
When the hyperplane $H$ in definition \ref{dfn:domebody} intersects the open ball bounded by $dev(\hat{\partial})$, the subset $D \cap dev^{-1}(H)$ is called the bottom of the dome body $D$.
\end{Definition}

\begin{Definition}[Bottom hyperplane]
The hyperplane $H$ in $\mathbb{A}^n$ containing the developing image of the bottom of the dome body $D$ is called the bottom hyperplane of $D$, and it is denoted by $\mathcal{H}(D)$.
\end{Definition}

\begin{Definition}[Bottom half space]
 The closed ( resp. open ) half space bounded by the bottom hyperplane of $D$ which also contains the interior of $dev(D)$ is called the closed ( resp. open ) bottom half space of $D$, and it is denoted by $\mathcal{U}(D)$ ( resp. $\stackrel{\circ}{\mathcal{U}}(D)$ ).
\end{Definition}

\begin{Definition} [Locally strictly convex] A subset $S$ of $\hat{M}$ is called locally strictly convex if for any $p \in S$, there exists a neighborhood $U$ of $p$ in $S$ such that $dev(U)$ is strictly convex in $\mathbb{A}^n$.
\end{Definition}

Remark. Our normalization of the boundary to a standard sphere in Theorem \ref{thm:convex_main} is to make sure that a neighborhood of the boundary is locally strictly convex.

The first key observation of dome bodies is the following.
\begin{Proposition} \label{prop:onebdry}
A dome body $D$ of $\hat{M}$ meets only one component of $\partial\hat{M}$, namely $\hat{\partial}$, which intersects $D$ on the top of $D$.
\end{Proposition}

\begin{Proof}
If $D$ meets another component $\hat{\partial}'$ of $\partial\hat{M}$ other than $\hat{\partial}$, we have the following contradiction:
\begin{itemize}
\item $\hat{\partial}'$ cannot meet the top of $D$, since $\hat{\partial}'$ is different from $\hat{\partial}$;
\item $\hat{\partial}'$ cannot meet the interior of $D$, since those are interior points of $\hat{M}$;
\item $\hat{\partial}'$ cannot meet the bottom of $D$, since any neighborhood of an intersection point would contain a small open half ball which makes $\hat{\partial}'$ fail to be locally strictly convex.
\end{itemize}
Therefore $D$ can only meet one component of $\partial\hat{M}$.
\end{Proof}

Remark. Proposition \ref{prop:onebdry} also says that all the points in the bottom of a dome body except those also in the top are interior points of $\hat{M}$. Therefore the developing map restricted to a dome body is a homeomorphism onto a closed solid semi ball minus some closed subset in the bottom hyperplane.

The second key observation is the following.
\begin{Proposition}\label{prop:domebody_compact}
Dome bodies are compact.
\end{Proposition}

Remark. Proposition \ref{prop:domebody_compact} says that the developing map restricted to a dome body is a homeomorphism onto a closed solid semi ball.

We postpone the proof to the next subsection. Assuming this, we have the following lemmas.

\begin{Lemma}\label{lemma:existence_of_max}
Let $\mathcal{D}ome$ be the set of dome bodies together with the partial order defined by inclusion. Then there exists a maximal element in $\mathcal{D}ome$.
\end{Lemma}

\begin{Proof}
Let $\mathcal{C}$ be a chain of dome bodies in $\mathcal{D}ome$. Then all the dome bodies in $\mathcal{C}$ share a common component $\hat{\partial}$ of $\partial\hat{M}$, and hence their developing images are all contained in the solid ball bounded by $dev(\hat{\partial})$. Therefore the set of their volumes $\{\, vol(dev(D)) \,|\, D \in \mathcal{C} \,\}$ is bounded above and let us call the supremum $V_{sup}$.

If $V_{sup}$ can be obtained by some $D$ in $\mathcal{C}$, then $D$ is obviously an upper bound for $\mathcal{C}$.

If $V_{sup}$ is not obtainable by any $D$ in $\mathcal{C}$, we can extract a sequence $(D_i)$ in $\mathcal{C}$ such that $vol(D_i)$ converges to $V_{sup}$. By passing to a subsequence, we can assume that the sequence of open half spaces $(\stackrel{\circ}{\mathcal{U}}(D_i))$ converges to some open half space $U_{\infty}$ in $\mathbb{A}^n$. Let $D_{\infty}$ be the closure of the union of $\{\stackrel{\circ}{D_i}\}$. Then the interior of $D_{\infty}$ is homeomorphic to the open solid semi ball which is the intersection of $U_{\infty}$ and the open solid ball bounded by $dev(\hat{\partial})$. Therefore $D_{\infty}$ is a dome body. Given any $D$ in $\mathcal{C}$, since $vol(D) < V_{sup}$, there exists some $D_i$ in the sequence we chose before such that $vol(D) < vol(D_i)$, therefore $D$ is contained in $D_i$ and hence is also contained in $D_{\infty}$. Then we have found an upper bound $D_{\infty}$ for $\mathcal{C}$ in this case as well.

Therefore, by Zorn's lemma there exists a maximal element in $\mathcal{D}ome$.
\end{Proof}

\begin{Lemma}\label{lemma:max_is_a_ball}
A maximal dome body $D_{max}$ is homeomorphic via $dev$ to the closed solid ball bounded by $dev(\hat{\partial})$, where $\hat{\partial}$ is the component of $\partial\hat{M}$ that meets $D_{max}$.
\end{Lemma}

\begin{Proof} Suppose $dev(D_{max})$ is not the whole closed solid ball $\mathcal{B}$ bounded by $dev(\hat{\partial})$. By Proposition \ref{prop:domebody_compact}, $dev(D_{max})$ is a closed solid semi ball. Let $V$ be a small open neighborhood of $D$. Since $dev$ is a local homeomorphism and that only the top of $D_{max}$ are boundary points of $\hat{M}$, $dev(V) \cap \mathcal{B}$ is an open subset in $\mathcal{B}$ that contains the closed solid semi ball $dev(D_{max})$. Therefore there exists a larger open solid semi ball containing $dev(D_{max})$ which is contained in $dev(V)$. Then the closure of the lift of it in $V$ defines a dome body that properly contains $D_{max}$. This contradicts to that $D_{max}$ is maximal. Therefore a maximal dome body $D_{max}$ is homeomorphic to the solid ball $\mathcal{B}$.
\end{Proof}

\begin{Lemma}\label{lemma:M_equals_D}
$\hat{M}$ is equal to a maximal dome body.
\end{Lemma}
\begin{Proof}
Since a maximal dome body $D_{max}$ is homeomorphic to the closed solid ball bounded by $dev(\hat{\partial})$, where $\hat{\partial}$ is the component of $\partial\hat{M}$ that meets $D_{max}$, it is straightforward to see that $D_{max}$ is both open and closed in $\hat{M}$, and hence a connected component of $\hat{M}$. Therefore $D_{max}$ is equal to $\hat{M}$.
\end{Proof}

Now we can finish the proof of Theorem \ref{thm:convex_main}.

\begin{Proof}[Proof of Theorem \ref{thm:convex_main}]
By Lemma \ref{lemma:M_equals_D} $\hat{M}$ is equal to a maximal dome body and hence $\partial\hat{M}$ has only one component by Propersition \ref{prop:onebdry}.
Therefore $\hat{M}$ is only a simple cover of $M$, i.e. $M$ is homeomorphic to $\hat{M}$ which is just a maximal dome body, which itself is homeomorphic to a closed solid ball under the developing map by Lemma \ref{lemma:max_is_a_ball}.
\end{Proof}

\subsection{Dome bodies are compact}
Here we are going to prove :\\
\textbf{Proposition \ref{prop:domebody_compact} }{\it Dome bodies are compact.}

\begin{Proof}
  For any $x,y \in \mathbb{A}^n$, we will use the notation $[x,y]$ to denote the straight line segment in $\mathbb{A}^n$ from $x$ to $y$ and use $[x, y)$ to denote the half open segment from $x$ to $y$.

  We will prove Proposition \ref{prop:domebody_compact} by contradiction in four steps.

  \noindent$\bullet$ Step 1, finding an incomplete geodesic.

  Suppose we have a dome body $D$ which is not compact. Let $C$ be the intersection of the top and the bottom of $D$, which is just a $n-2$ sphere. Then there exists two distinct points $a$ and $b$ in $C$ such that not the whole segment $[dev(a),dev(b)]$ can be lifted to a segment in $\hat{M}$ starting at $a$. For otherwise if any segment with endpoints in $dev(C)$ can be lifted to $\hat{M}$, the closure of $dev(D)$ can be lifted, and hence $dev(D)$ is closed and $D$ is homeomorphic to a closed solid semi ball, which contradicts to that $D$ is not compact.

  \begin{figure}[h]
  \begin{center}
  \scalebox{0.6} {\epsffile{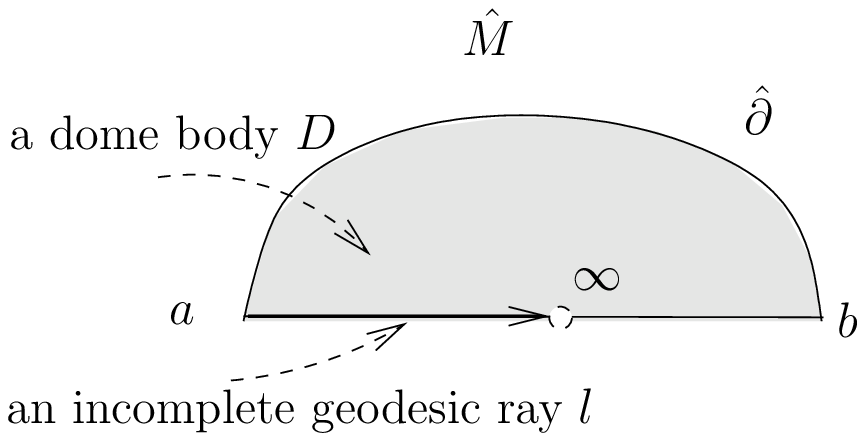}}
  \scalebox{0.35} {\epsffile{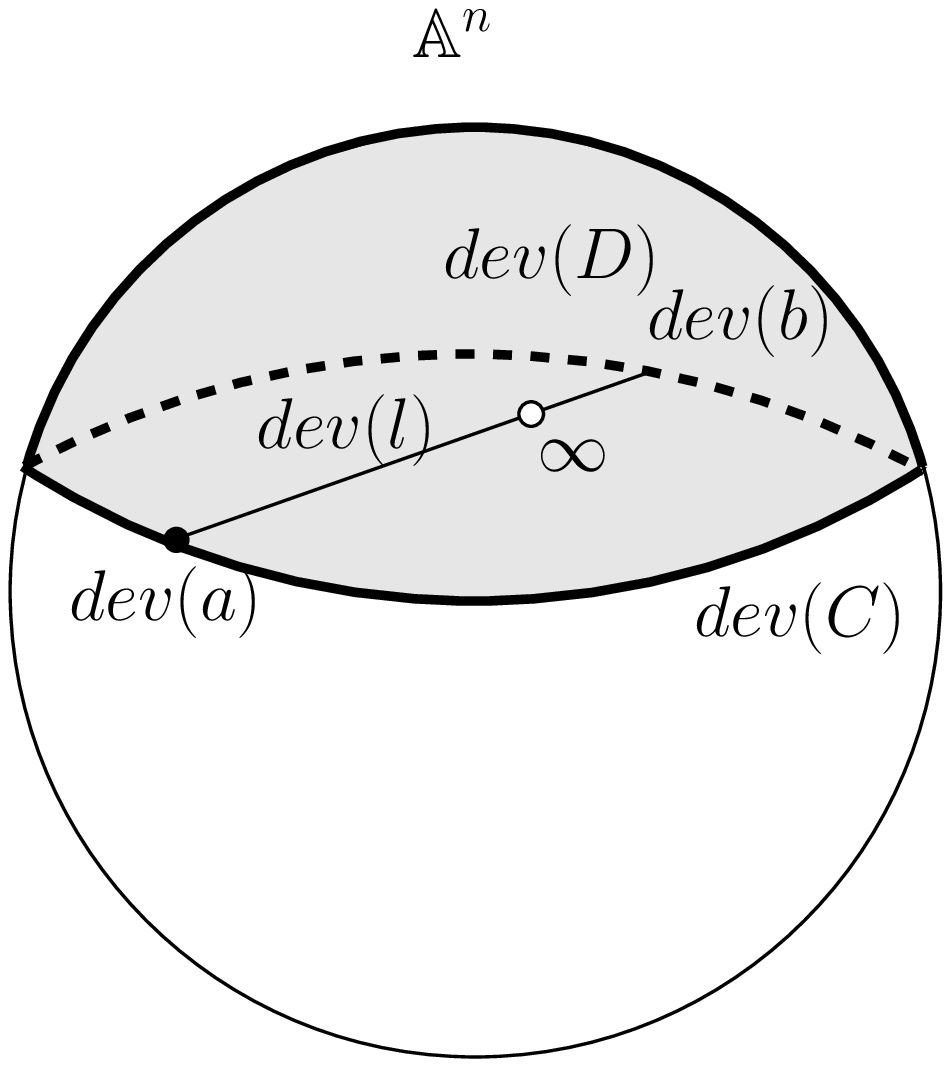}}
  \end{center}
  \caption{An incomplete geodesic ray in the bottom of a noncompact dome body.} \label{fig:inc_geodesic}
  \end{figure}

  Then there exists a maximal half open subsegment $[dev(a), \infty)$ in $[dev(a), dev(b)]$ such that it can be lifted, and we use $l$ to denote its lift starting at $a$. Note that $\infty$ cannot be $dev(b)$, so $\infty$ lives in the interior of $[dev(a), dev(b)]$. Then $l$ necessarily leaves any compact set of $\hat{M}$. We call $l$ an incomplete geodesic ray, since $dev(l)$ is a bounded straight ray in $\mathbb{A}^n$ and $l$ leaves any compact set of $\hat{M}$ ( See Figure \ref{fig:inc_geodesic} ).

  \noindent$\bullet$ Step 2, finding a sequence of dome bodies.

  Since $M$ is compact, the projection $\bar{l}$ of $l$ to $M$ is a recurrent geodesic, i.e. it will not stay in any small compact set eventually. Let $\bar{p}$ be an accumulation point of $\bar{l}$ in $M$. Since the boundary is locally strictly convex, $\bar{p}$ is necessarily an interior point of $M$. For otherwise if $\bar{p}$ is on $\partial M$, we can take a small compact semi solid ball neighborhood of $\bar{p}$, but then after $\bar{l}$ enters this compact neighborhood $\bar{l}$ will stay in it eventually, which is a contradiction.

  \begin{figure}[h]
  \begin{center}
  \scalebox{.8} {\epsffile{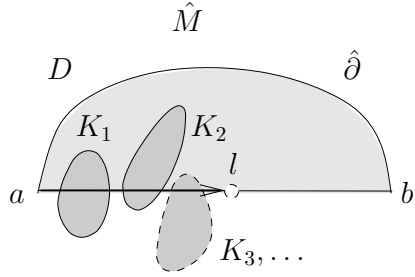}}
  \end{center}
  \caption{The sequence of lifts of $K$.} \label{fig:Ks}
  \end{figure}

  Then we can choose a small compact convex neighborhood $\bar{K}$ of $\bar{p}$ ( By convex, we mean the developing image of a lift of $\bar{K}$ in $\mathbb{A}^n$ is convex ), which is also contained in the interior of $M$, such that $\bar{l}$ enters and leaves $\bar{K}$ infinitely many times. Therefore $l$ meets a sequence of lifts $(K_i)$ of $\bar{K}$ in the universal cover $\hat{M}$ ( See Figure \ref{fig:Ks} ).

  Note that $\bar{K}$ does not contain any boundary point of $\hat{M}$, so $K_i \cap D$ does not meet the top of $D$. The bottom hyperplane $\mathcal{H}(D)$ separates the interior of $dev(K_i)$ into two open semi balls. By the convexity of $dev(K_i)$ and the interior of $dev(D)$, the inverse image of one of these open semi balls in $K_i$ is completely contained in $D$ while the inverse image of the other one in $K_i$ is disjoint from $D$.

  Since $\bar{K}$ was chosen to be a small compact convex neighborhood and $\{K_i\}$ are different lifts of $\bar{K}$, they are disjoint. Since $\bar{l}$ enters and leaves $\bar{K}$ infinitely many times, $l$ enters and leaves the sequence $(K_i)$ successively. Since $dev(l)$ is bounded in $\mathbb{A}^n$, the length of $dev(K_i \cap l)$ must go to $0$. Therefore we have

  \begin{Lemma} \label{lemma:Kcaplsmall}
  The length of $dev(l \cap K_i)$ goes to $0$ as $i$ goes to infinity.
  \end{Lemma}

  \begin{figure}[h]
  \begin{center}
  \scalebox{.8} {\epsffile{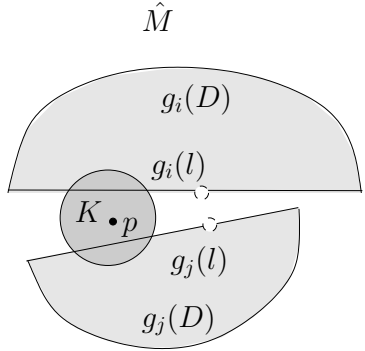}}
  \end{center}
  \caption{The sequence of $\{g_i(D)\}$ meeting $K$.} \label{fig:Ds}
  \end{figure}

  Now we fix a lift $K$ of $\bar{K}$ in $\hat{M}$. By using the deck transformation $g_i$ which takes $K_i$ to $K$, we have a sequence of dome bodies $\left(g_i(D)\right)$, all of which meet $K$ ( See Figure \ref{fig:Ds} ). Then $g_i(D) \cap K$ does not meet the top of $g_i(D)$. Once again the bottom hyperplane $\mathcal{H}(g_i(D))$ separates the interior of $dev(K)$ into two open semi balls. By the convexity of $dev(K)$ and the interior of $dev(g_i(D))$, the inverse image of one of these open semi balls in $K$ is completely contained in $g_i(D)$ while the inverse image of the other one in $K$ is disjoint from $g_i(D)$. Note that we also have a lift $p$ of $\bar{p}$ in $K$ which is an accumulation point of $\left(g_i(l)\right)$.

  \noindent$\bullet$ Step 3, stacking a ``small'' dome body and a ``large'' dome body.

  We will use $\mathcal{L}_i$ to denote the line in $\mathbb{A}^n$ containing $dev(g_i(l))$ in the following.

  Let us look at the developing images $\left\{dev(g_i(D))\right\}$ in $\mathbb{A}^n$, the closure of which are just closed solid semi balls. Since $dev(K)$ is compact and $dev(g_i(D))$ intersects $dev(K)$, by passing to a subsequence we can assume that the sequence of lines $(\mathcal{L}_i)$ converges to a line $l_{\infty}$ passing through $dev(p)$, and by passing to another subsequence we can also assume that the sequence of closed bottom half spaces $ \left( \mathcal{U}(g_i(D)) \right)$ converges to some closed half space $U_{\infty}$ of $\mathbb{A}^n$, the boundary hyperplane $H_{\infty}$ of which necessarily contains $dev(p)$ and $l_{\infty}$.

  \begin{figure}[h]
  \begin{center}
  \scalebox{0.8} {\epsffile{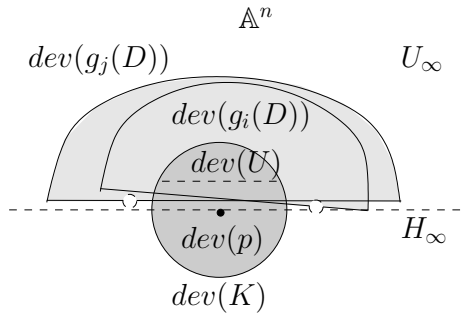}}
  \end{center}
  \caption{A subsequence of $\{dev(g_i(D))\}$ containing $dev(U)$.} \label{fig:devDs}
  \end{figure}

  If we fix a hyperplane which is parallel and close to $H_{\infty}$ and is contained in the interior of $U_{\infty}$, then it still intersects the interior of $dev(K)$ and it separates the interior of $dev(K)$ into two open solid semi balls, one of which is contained in $U_\infty$ and its lift in $K$ is denoted by $U$. We call $U$ the upper part of $K$.

  Since the sequence of closed bottom half spaces $ \left( \mathcal{U}(g_i(D)) \right)$ converges to $U_{\infty}$, when $i$ is large enough every $dev(g_i(D))$ must contain $dev(U)$ and hence $g_i(D)$ contains $U$ ( See Figure \ref{fig:devDs}, $dev(g_i(D))$ contains the upper part $dev(U)$ of $dev(K)$ bounded by the dash plane parallel to $H_{\infty}$ in the figure ).

  So by passing to a subsequence again, we can assume
  \begin{Lemma}\label{lemma:containU}
  Any dome body in the sequence $(g_i(D))$ contains the upper part $U$ of $K$.
  \end{Lemma}

  Note that affine transformations preserve the ratio of the length of two parallel segments. If we use $|.|$ to denote the length, then
  \begin{eqnarray*}
  \frac{\left|\, dev(l \cap K_i) \,\right|}{\left|\, [dev(a), dev(b)] \,\right|} & = & \frac{\left|\, dev\left(g_i(l \cap K_i)\right) \,\right|}{\left|\, [dev(g_i(a)), dev(g_i(b))] \,\right|} \\
  & = & \frac{\left|\, dev\left(g_i(l)\cap K\right) \,\right|}{\left|\, [dev(g_i(a)), dev(g_i(b))] \,\right|}.
  \end{eqnarray*}

  On one hand, by Lemma \ref{lemma:Kcaplsmall} the lhs goes to $0$ as $i$ goes to infinity; on the other hand, when $i$ is large enough, $g_i(l)$ will be close to $l_{\infty}$ which passes through the point $p$ in the interior of $K$, hence the numerator on the rhs $$\left| \, dev(g_i(l) \cap K) \,\right| \quad \to \quad \left| \, dev(l_{\infty} \cap K) \,\right| > 0,$$ and hence the denominator $$\left| \, [dev(g_i(a)), dev(g_i(b))] \,\right| \quad \to \quad \infty$$ as $i$ goes to infinity.

  Similarly, let $c_i$ and $d_i$ be the endpoints of $l \cap K_i$ and replace the pair of points $a$ and $b$ above by the pair of $a$ and $d_i$ ( resp. $c_i$ and $b$ ). Since $dev(c_i)$ ( resp. $dev(d_i)$ ) is close to the point $\infty$ in the interior of $[dev(a), dev(b)]$, we actually have
  \begin{Lemma}\label{lemma:abgoestoinfty}
  $dev(g_i(a))$ ( resp. $dev(g_i(b))$ ) goes to infinity in $\mathbb{A}^n$ as $i$ goes to infinity.
  \end{Lemma}

  From now on we will fix $i$ and take $j$ arbitrarily large. By ``small'', we mean $dev(g_i(D))$ is bounded; by ``large'', we mean both $dev(g_j(a))$ and $dev(g_j(b))$ can be arbitrarily far away from $dev(g_i(D))$; by ``stacking'', we mean they both contain a common part $U$ of $K$.

  \noindent$\bullet$ Step 4, finding a contradiction.

  \begin{figure}[h]
  \begin{center}
  \scalebox{.8} {\epsffile{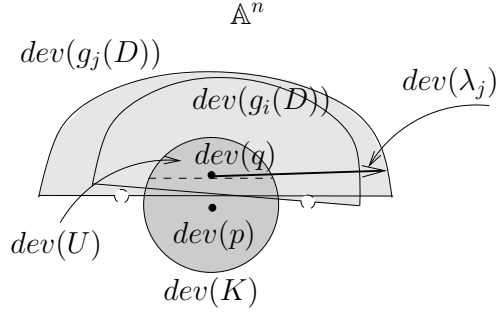}}
  \end{center}
  \caption{A long segment $\lambda_j$ in $g_j(D)$.} \label{fig:lambdaj}
  \end{figure}

  Now let us fix a point $q$ in the interior of $U$. Starting from $dev(q)$ there are two geodesic segments parallel to the line $\mathcal{L}_j$, each of which has the other endpoint on the developing image of the top of $g_j(D)$. Note that they point in opposite directions, so there is at least one not pointing towards the bottom hyperplane of $g_i(D)$, and let us use $\lambda_j$ to denote its lift in $g_j(D)$ ( See Figure \ref{fig:lambdaj} ).

  We are going to prove that $\left| dev(\lambda_j) \right|$ can be arbitrarily large if $j$ is large enough. For this, without loss of generality, let us assume that $dev(\lambda_j)$ points in the direction from $dev(g_j(a))$ to $dev(g_j(b))$. Let us take the line through $dev(q)$ perpendicular to the line $\mathcal{L}_j$ containing $dev(g_j(l))$ with the intersection point $\nu_j$ ( note that $\nu_j$ may not be lifted to $\hat{M}$ ). By the convexity of both $dev(K)$ and the interior of $dev(g_j(D))$, the inverse image of this line intersects $K$ at a point $m_j$ such that $dev(q)$ is in the interior of the segment $[dev(m_j), \nu_j]$ and hence $m_j$ is contained in $g_j(D)$ ( See Figure \ref{fig:triangle} ).

  \begin{figure}[h]
  \begin{center}
  \scalebox{0.4} {\epsffile{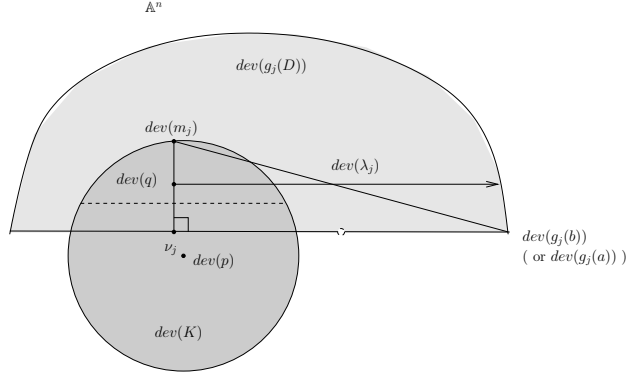}}
  \end{center}
  \caption{Why is $\lambda_j$ long?} \label{fig:triangle}
  \end{figure}

  By the convexity of the interior of $dev(g_j(D))$ again, the segment\\ $[dev(m_j), dev(g_j(b))]$ is contained in $dev(g_j(D))$. Using similar triangles, we have
  $$ \left| \, dev(\lambda_j) \,\right| \quad \geq \quad \frac{\left|\, [dev(m_j),dev(q)] \,\right|}{\left|\, [dev(m_j),\nu_j] \,\right|} \, \cdot \, \left|\, [dev(g_j(b)),\nu_j] \,\right|. $$

  Since $(\mathcal{L}_j)$ converges to $l_{\infty}$ and $$\frac{\left|\, [dev(m_j),dev(q)] \,\right|}{\left|\, [dev(m_j),\nu_j] \,\right|}$$ depends continuously on them, it converges to some limit bounded away from $0$. Since $\nu_j$ is within bounded distance from $dev(K)$, by Lemma \ref{lemma:abgoestoinfty} $$\left|\, [dev(g_j(b)),\nu_j] \,\right| \quad \to \quad \infty.$$ Therefore we have

  \begin{Lemma}\label{lemma:lambdajlong}
  $\left|\, dev(\lambda_j) \,\right| \to \infty$ as $j$ goes to infinity.
  \end{Lemma}

  Now $s \doteq \lambda_j \cap g_i(D)$ is a nonempty segment in $g_i(D)$, since both $\lambda_j$ and $g_i(D)$ are convex. Let $r$ be the other endpoint of $s$ opposite to $q$.

  On one hand, since $[dev(q), dev(r)]$ is contained in $dev(g_i(D))$ and hence is bounded. If we choose $j$ large enough in the first place, we have $$\left|\, dev(\lambda_j) \,\right| \quad > \quad \left|\, [dev(q), dev(r)] \,\right|. $$ Then $r$ must live in the interior of $\lambda_j$ and hence $r$ is an interior point of $\hat{M}$ if you look at it in $g_j(D)$; on the other hand, by the choice of $\lambda_j$, $r$ is on the top of $g_i(D)$, and hence $r$ is a boundary point of $\hat{M}$ if you look at it in $g_i(D)$. Therefore we have a contradiction.
\end{Proof}

\subsection{Extension to projective case}
In this subsection, we are going to extend our result to the projective case and prove\\
\textbf{Theorem \ref{thm:proj} }
{\it For $n \geq 3$, let $(M, \partial)$ be a compact projective n-manifold with boundary $\partial$ homeomorphic to $S^{n-1}$. If
\begin{itemize}
  \item $dev$ restricted to some lift $\hat{\partial}$ of $\partial$ is an embedding,
  \item $dev(\hat{\partial})$ is contained in an affine patch $\mathbb{A}^n$,
  \item $dev$ maps a neighborhood of $\hat{\partial}$ to the closure of the bounded part of $\mathbb{A}^n \backslash dev(\hat{\partial})$,
\end{itemize}
then $(M, \partial)$ is homeomorphic to $(D^n, S^{n-1})$. }

\begin{Proof}
All the arguments in \S4.1 work through in this case as long as we can establish the projective version of Proposition \ref{prop:domebody_compact}: Dome bodies are compact.

Once again we will prove it by contradiction. If there is a noncompact dome body $D$, we can find an incomplete geodesic ray $l$ in the bottom of $D$, which goes through a sequence of $(K_i)$, where $K_1, K_2, \dots$ are lifts of some compact convex neighborhood of an accumulation point ( See Figure \ref{fig:projKs} ).

\begin{figure}[h]
\begin{center}
\scalebox{.8} {\epsffile{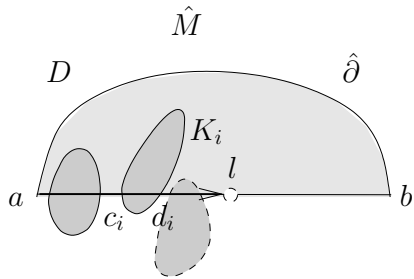}}
\end{center}
\caption{The sequence of lifts of $K$ in the projective case.} \label{fig:projKs}
\end{figure}

Let $c_i$ and $d_i$ be the endpoints of $l \cap K_i$. For the four points $a, c_i, d_i, b$ in a line, we can consider the cross-ratio
$$
(a,c_i;d_i,b) \doteq \frac{\left| \, [dev(a),dev(d_i)] \, \right| \cdot \left| \, [dev(c_i), dev(b)] \, \right|}{\left| \, [dev(a),dev(b)] \, \right| \cdot \left| \, [dev(c_i), dev(d_i)] \, \right|}
$$
which is a projective invariant.

Since $\left| \, [dev(a),dev(d_i)] \, \right|$, $\left| \, [dev(c_i),dev(b)] \, \right|$ and $\left| \, [dev(a),dev(b)] \, \right|$ are bounded above and below away from $0$ uniformly for all $i$ and $\left| \, [dev(c_i), dev(d_i)] \, \right| \to 0$ as $i$ goes to infinity ( See Figure \ref{fig:projKs} ), we have
\begin{Lemma}\label{lemma:crgoestoinf}
$(a,c_i;d_i,b) \to \infty$ as $i$ goes to infinity.
\end{Lemma}

If we fix a $K$ and use the deck transformations $(g_i)$, where $g_i$ takes $K_i$ to $K$, then $K$ intersects the sequence of dome bodies $(g_i(D))$ again. From Lemma \ref{lemma:crgoestoinf}, we have
\begin{Lemma}\label{lemma:crgoestoinf1}
$(g_i(a),g_i(c_i);g_i(d_i),g_i(b)) \to \infty$ as $i$ goes to infinity.
\end{Lemma}

The difference between the projective case and the affine case is that we might not be able to see the whole $dev(g_i(D))$ in an affine patch containing $dev(K)$. This is not a big issue. We can still start with an affine patch containing $dev(K)$. By passing to a subsequence we can still assume that the sequence $(\mathcal{L}_i)$, where $\mathcal{L}_i$ is the line containing $dev(g_i(l))$, converges to some line $l_{\infty}$, $\left(\mathcal{U}(g_i(D)) \right)$ converges to some closed half space $U_{\infty}$ and $\left(\mathcal{H}(g_i(D))\right)$ converges to some hyperplane $H_{\infty}$. Hence by passing to a subsequence we can still have the projective version of Lemma \ref{lemma:containU}: Any dome body in the sequence $(g_i(D))$ contains an upper part $U$ of $K$.

Once again, from now on we will fix $i$ and a point $q$ in the interior of $U$. Starting from $dev(q)$ there are two geodesic segments parallel to the line through $dev(g_j(l))$ with the other endpoint on the developing image of the top of $g_j(D)$ ( We can switch to an affine patch containing the closure of $dev(g_j(D))$ to do this ). Note that they point in opposite directions, so there is at least one not pointing towards the bottom hyperplane of $g_i(D)$, and let us use $\lambda_j$ to denote its lift in $g_j(D)$.

Once again we consider the subsegment $s \doteq \lambda_j \cap g_i(D)$. If $s$ is properly contained in $\lambda_j$, we will have the same contradiction as before: the other endpoint $r$ of $s$ opposite to $q$ is both an interior point and a boundary point of $\hat{M}$. Therefore $s$ must equal to $\lambda_j$, which means that $g_i(D)$ and $g_j(D)$ share a common component of $\partial\hat{M}$.

If we switch to an affine patch containing the closed solid ball bounded by the developing image of the component of $\partial\hat{M}$ that meets $g_i(D)$. Then $\left\{\, dev(g_j(D)) \,\right\}$ are all contained in the same solid ball. Since the lengths \\ $\left| \, [dev(g_j(a)),dev(d_j)] \, \right|$ and $\left| \, [dev(g_j(c_j)),dev(g_j(b))] \, \right|$ are uniformly bounded above while $\left| \, [dev(g_j(a)),dev(g_j(b))] \, \right|$ and $\left| \, [dev(g_j(c_j)),dev(g_j(d_j))] \, \right|$ are uniformly bounded away from $0$ for all $j$, and hence
\begin{eqnarray*}
 & & (g_j(a),g_j(c_i);g_j(d_i),g_j(b))\\
 & =  & \frac{\left| \, [dev(g_j(a)),dev(g_j(d_j))] \, \right| \cdot \left| \, [dev(g_j(c_j)), dev(g_j(b))] \, \right|}{\left| \, [dev(g_j(a)),dev(g_j(b))] \, \right| \cdot \left| \, [dev(g_j(c_j)), dev(g_j(d_j))] \, \right|}
\end{eqnarray*}
is uniformly bounded above. This contradicts to Lemma \ref{lemma:crgoestoinf1}.
\end{Proof}